\documentclass[11pt, a4paper]{article}

\usepackage{array}
\usepackage{multirow}

\makeatletter
\newcommand{\thickhline}{%
    \noalign {\ifnum 0=`}\fi \hrule height 1pt
    \futurelet \reserved@a \@xhline
}
\newcolumntype{"}{@{\hskip\tabcolsep\vrule width 1pt\hskip\tabcolsep}}
\makeatother

\usepackage{graphics}
\usepackage{graphicx}
\usepackage{epsfig}
\usepackage{subfigure}
\usepackage{amssymb}
\usepackage{amsthm}
\usepackage{mathrsfs}
\usepackage{hhline}
\usepackage{longtable}
\usepackage{amsmath}
\usepackage{array}
\usepackage{float}
\usepackage{multirow}
\usepackage{cite}
\usepackage{enumerate}
\usepackage{xcolor}

\def\ms{\medskip}
\def\nt{\noindent}

\definecolor{vividviolet}{rgb}{0.62, 0.0, 1.0}

\newtheoremstyle{de}
  {10pt}          
  {10pt}  
  {\rm}  
  {}
  {\bf}  
  {. }    
  { }    
  {}     
\theoremstyle{de}

\newtheorem{de}{Definition}[section]

\newtheorem{rem}[de]{Remark}

\newtheorem{question}{Question}[section]
\newtheorem{definition}{Definition}[section]
\newtheoremstyle{theorem}
  {10pt}          
  {10pt}  
  {\it}  
  {}
  {\bf}  
  {. }    
  { }    
  {}     
\theoremstyle{theorem}
\newtheorem{theorem}{Theorem}[section]
\newtheorem{lemma}[theorem]{Lemma}

\newtheorem{corollary}[theorem]{Corollary}

\numberwithin{equation}{section}

\graphicspath{%
    {converted_graphics/}
    {/}
}

\begingroup\makeatletter\ifx\SetFigFont\undefined%
\gdef\SetFigFont#1#2#3#4#5{%
  \reset@font\fontsize{#1}{#2pt}%
  \fontfamily{#3}\fontseries{#4}\fontshape{#5}%
  \selectfont}%
\fi\endgroup%

\begin{document}
\begin{center}
{\mathversion{bold}\Large \bf On join product and local antimagic chromatic number of regular graphs}

\bigskip
{\large  Gee-Choon Lau$^{a,}$\footnote{Corresponding author.}, Wai Chee Shiu$^{b}$}\\

\medskip

\emph{{$^a$}Faculty of Computer \& Mathematical Sciences,}\\
\emph{Universiti Teknologi MARA (Johor Branch, Segamat Campus),}\\
\emph{85000, Malaysia.}\\
\emph{geeclau@yahoo.com}\\

\medskip
\emph{{$^b$}Department of Mathematics,}\\
\emph{The Chinese University of Hong Kong,}\\
\emph{Shatin, Hong Kong, China.}\\
\emph{wcshiu@associate.hkbu.edu.hk}\\

\end{center}

\begin{abstract} Let $G = (V,E)$ be a connected simple graph of order $p$ and size $q$. A graph $G$ is called local antimagic if $G$ admits a local antimagic labeling. A bijection $f : E \to \{1,2,\ldots,q\}$ is called a local antimagic labeling of $G$ if for any two adjacent vertices $u$ and $v$, we have $f^+(u) \ne f^+(v)$, where $f^+(u) = \sum_{e\in E(u)} f(e)$, and $E(u)$ is the set of edges incident to $u$. Thus, any local antimagic labeling induces a proper vertex coloring of $G$ if vertex $v$ is assigned the color $f^+(v)$. The local antimagic chromatic number, denoted $\chi_{la}(G)$, is the minimum number of induced colors taken over local antimagic labeling of $G$. Let $G$ and $H$ be two vertex disjoint graphs. The join graph of $G$ and $H$, denoted $G \vee H$, is the graph with $V(G\vee H) = V(G) \cup V(H)$ and $E(G\vee H) = E(G) \cup E(H) \cup \{uv \,|\, u\in V(G), v \in V(H)\}$. In this paper, we show the existence of non-complete regular graphs with arbitrarily large order, regularity and local antimagic chromatic numbers.
 \\


\end{abstract}

\section{Introduction}

Let $G = (V,E)$ be a connected simple graph of order $p$ and size $q$. A bijection $f : E \to \{1,2,\ldots,q\}$ is called a {\it local antimagic labeling} of $G$ if for any two adjacent vertices $u$ and $v$, we have $f^+(u) \ne f^+(v)$, where $f^+(u) = \sum_{e\in E(u)} f(e)$, and $E(u)$ is the set of edges incident to $u$. Thus, any local antimagic labeling induces a proper vertex coloring of $G$ if vertex $v$ is assigned the color $f^+(v)$. If $f$ induces $t$ distinct colors, we say $f$ is a {\it local antimagic $t$-coloring} of $G$. The {\it local antimagic chromatic number} of $G$, denoted $\chi_{la}(G)$, is the minimum number of induced colors taken over local antimagic labelings of $G$ \cite{Arumugam}. 

\ms\nt Let $G$ and $H$ be two vertex disjoint graphs. The {\it join graph} of $G$ and $H$, denoted $G \vee H$, is the graph $V(G\vee H) = V(G) \cup V(H)$ and $E(G\vee H) = E(G) \cup E(H) \cup \{uv \,|\, u\in V(G), v \in V(H)\}$. In~\cite{LNS, LNS-dmgt, LauShiuNg, YBY}, the local antimagic chromatic number of the join of graphs are determined. However, very few results on regular graphs are obtained. Motivated by this, we investigate the local antimagic chromagic number of the join of regular graphs. In Section 2, we completely determined the local antimagic chromatic number of regular graphs of order at most 8. In Section 3, we constructed and determined the local antimagic chromatic number of infinitely many regular graphs using join product operation repeatedly. Consequently, we showed the existence of non-complete regular graphs with arbitrarily large order, regularity and local antimagic chromatic numbers.

\ms\nt For convenient, we shall use $[a,b]$ to denote the set of integers from $a$ to $b$, where $a\le b$ are integers. We first give some useful lemmas and definitions.

\begin{lemma}[\hspace*{-1mm}\cite{LNS-dmgt}]\label{lem-2part}
Let $G$ be a graph of size $q$. Suppose there is a local antimagic labeling of $G$ inducing a $2$-coloring of $G$ with colors $x$ and $y$, where $x<y$.  Let $X$ and $Y$ be the sets of vertices colored $x$ and $y$, respectively, then $G$ is a bipartite graph with bipartition $(X,Y)$ and $|X|>|Y|$. Moreover,
\begin{equation*} x|X|=y|Y|= \frac{q(q+1)}{2}.\end{equation*}
\end{lemma}

\nt Since bipartition of a connected bipartite graph is unique, Lemma~\ref{lem-2part} implies that

\begin{corollary}\label{cor-2part} Suppose $G$ is a connected bipartite graph of $q$ edges with bipartition $(V_1,V_2)$. If $\chi_{la}(G)=2$, then $|V_1|\ne |V_2|$ and $\binom{q+1}{2}$ is divisible by both $|V_1|$ and $|V_2|$.\end{corollary}

\nt In what follows, let $G-e$ denote the graph $G$ with an edge $e$ deleted.

\begin{lemma}[\hspace*{-1mm}\cite{LNS-dmgt}]\label{lem-reg}  Suppose $G$ is a regular graph of size $q$. If $f$ is a local antimagic labeling of $G$, then $g = q + 1 - f$ is also a local antimagic labeling of $G$ with $c(f)= c(g)$.  Moreover, suppose $c(f)=\chi_{la}(G)$ and if $f(e)=1$ or $f(e)=q$, then $\chi_{la}(G-e)\le \chi_{la}(G)$. \end{lemma}

\nt A {\it labeling matrix} $M$ for an edge labeling $f$ of $G$ is a symmetric matrix whose rows and columns are named by the vertices of $G$ and the $(u,v)$-entry is $f(uv)$ if $uv \in E$, and is $*$ otherwise. Sometimes, we call this matrix a {\it labeling matrix of $G$}. In other words, suppose $A$ is an adjacency matrix of $G$ and $f$ is a labeling of $G$, then a labeling matrix for $f$ is obtained from $A=(a_{u,v})$ by replacing $a_{u,v}$ by $f(uv)$ if $a_{u,v}=1$ and by $*$ if $a_{u,v}=0$. This concept was first introduced by Shiu, {\it et al.} in \cite{Shiu1998}.

\ms\nt Let $G$ be a graph with vertex list $\{u_1, u_2,\dots, u_p\}$. Suppose $f$ is a local antimagic labeling of $G$ and $M$ is the corresponding labeling matrix. Then the $i$-th row sum denoted by $r_i(M)$ (also the $i$-th column sum $c_i(M)$) of $M$ is $f^+(u_i)$. Here the label $*$ is treated as 0. Thus if we want to construct a local antimagic labeling of a graph $G$ of order $p$ and size $q$, then it is equivalent to
\begin{enumerate}[(1)]
\item replace all 1's lying in the upper triangular part of the adjacent matrix of $G$ by labels in $[1, q]$ bijectively;
\item replace all 0's lying in the upper triangular part or the diagonal of the adjacent matrix of $G$ by $*$'s;
\item fill the lower triangular part so that the matrix becomes symmetric; \end{enumerate}
such that the $i$-th row sum is not equal to the $j$-th row sum if $u_iu_j$ is an edge in $G$.

\section{Connected Regular Graphs of Order at Most 8}\label{sec-regular}

For a fixed $n$ with $3\le n\le 8$, there is only one connected 2-regular graph which is $C_n$. The local antimagic chromatic number were determined. We shall consider circulant graph $C_n(1, a_1, \dots, a_t)$ that is obtained from an $n$-cycle by joining each pair of vertices with distance belongs to $\{a_1,\dots, a_t\}$ in $C_n$, where $1<a_1<\cdots <a_t\le \lfloor n/2\rfloor$ be a strictly increasing sequence of length $t$. Note that when $t=0$, the sequence does not exist and we get $C_n(1) = C_n$.

\subsection{Connected Cubic Graphs of Order at Most 8}

In this subsection, we consider connected cubic graphs or their related graphs of order at most 8.
Before the consideration we recall some definitions of famous classes of graphs. The Cartesian product $C_n\times P_2$ is called a {\it prism}. A {\it m\"obius ladder} $M_{2n}$ of order $2n$ is obtained from a regular $2n$-gon by adding an edge between each pair of antipodal vertices. The regular $2n$-gon is called the {\it outer cycle} of $M_{2n}$.

\ms\nt Let $G$ be a cubic graph of order $2n$. When $n=2$, then $G\cong K_4$. When $n=3$, then $G\cong K_{3,3}$ or $C_3\times P_2$. When $n=4$, then $G$ is isomorphic to one of the following graphs.

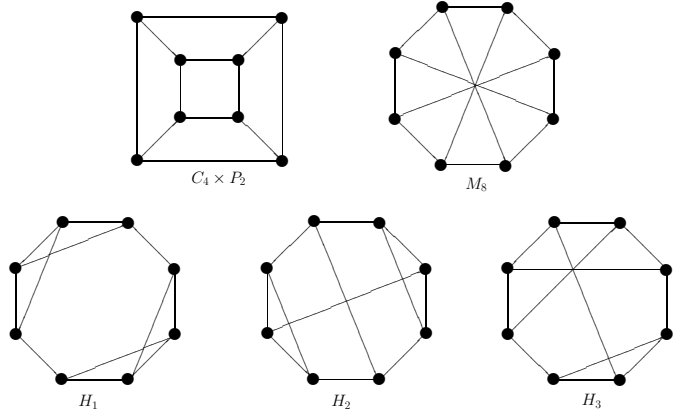
\begin{figure}[H]
\begin{center}
\scalebox{0.3}{
\setlength{\unitlength}{3947sp}%
\begingroup\makeatletter\ifx\SetFigFont\undefined%
\gdef\SetFigFont#1#2#3#4#5{%
  \reset@font\fontsize{#1}{#2pt}%
  \fontfamily{#3}\fontseries{#4}\fontshape{#5}%
  \selectfont}%
\fi\endgroup%
\begin{picture}(13815,8590)(2826,-9552)
\thicklines
{\color[rgb]{0,0,0}\put(6306,-6602){\line(-2,-5){930}}
}%
{\color[rgb]{0,0,0}\put(13174,-4411){\line( 1, 1){965}}
\put(14143,-3450){\line( 0, 1){1364}}
\put(14149,-2086){\line(-1, 1){965}}
\put(13188,-1117){\line(-1, 0){1364}}
\put(11824,-1111){\line(-1,-1){965}}
\put(10855,-2072){\line( 0,-1){1364}}
\put(10849,-3436){\line( 1,-1){965}}
\put(11810,-4405){\line( 1, 0){1364}}
}%
{\color[rgb]{0,0,0}\put(11843,-1100){\line( 2,-5){1330.069}}
}%
{\color[rgb]{0,0,0}\put(13171,-1110){\line(-2,-5){1321.793}}
}%
{\color[rgb]{0,0,0}\put(15506,-8917){\line( 1, 1){965}}
\put(16475,-7956){\line( 0, 1){1364}}
\put(16481,-6592){\line(-1, 1){965}}
\put(15520,-5623){\line(-1, 0){1364}}
\put(14156,-5617){\line(-1,-1){965}}
\put(13187,-6578){\line( 0,-1){1364}}
\put(13181,-7942){\line( 1,-1){965}}
\put(14142,-8911){\line( 1, 0){1364}}
}%
{\color[rgb]{0,0,0}\put(5510,-4328){\framebox(3000,3000){}}
}%
{\color[rgb]{0,0,0}\put(6410,-3428){\framebox(1200,1200){}}
}%
{\color[rgb]{0,0,0}\put(5510,-1328){\line( 1,-1){900}}
}%
{\color[rgb]{0,0,0}\put(5320,-8910){\line( 1, 1){965}}
\put(6289,-7949){\line( 0, 1){1364}}
\put(6295,-6585){\line(-1, 1){965}}
\put(5334,-5616){\line(-1, 0){1364}}
\put(3970,-5610){\line(-1,-1){965}}
\put(3001,-6571){\line( 0,-1){1364}}
\put(2995,-7935){\line( 1,-1){965}}
\put(3956,-8904){\line( 1, 0){1364}}
}%
{\color[rgb]{0,0,0}\put(8510,-1328){\line(-1,-1){900}}
}%
{\color[rgb]{0,0,0}\put(7610,-3428){\line( 1,-1){900}}
}%
{\color[rgb]{0,0,0}\put(5510,-4328){\line( 1, 1){900}}
}%
{\color[rgb]{0,0,0}\put(10860,-2073){\line( 5,-2){3297.586}}
}%
{\color[rgb]{0,0,0}\put(9171,-5582){\line( 2,-5){1325.172}}
}%
{\color[rgb]{0,0,0}\put(10536,-5612){\line( 2,-5){961.034}}
}%
{\color[rgb]{0,0,0}\put(8181,-6557){\line( 2,-5){950.690}}
}%
{\color[rgb]{0,0,0}\put(8196,-7922){\line( 5, 2){3292.241}}
}%
{\color[rgb]{0,0,0}\put(10526,-8902){\line( 1, 1){965}}
\put(11495,-7941){\line( 0, 1){1364}}
\put(11501,-6577){\line(-1, 1){965}}
\put(10540,-5608){\line(-1, 0){1364}}
\put(9176,-5602){\line(-1,-1){965}}
\put(8207,-6563){\line( 0,-1){1364}}
\put(8201,-7927){\line( 1,-1){965}}
\put(9162,-8896){\line( 1, 0){1364}}
}%
{\color[rgb]{0,0,0}\put(14136,-5597){\line( 2,-5){1325.172}}
}%
{\color[rgb]{0,0,0}\put(14151,-8912){\line( 5, 2){2353.448}}
}%
{\color[rgb]{0,0,0}\put(13191,-6602){\line( 1, 0){3300}}
}%
{\color[rgb]{0,0,0}\put(13176,-7967){\line( 1, 1){2332.500}}
}%
{\color[rgb]{0,0,0}\put(3966,-5612){\line(-2,-5){944.483}}
}%
{\color[rgb]{0,0,0}\put(3006,-6572){\line( 5, 2){2343.103}}
}%
{\color[rgb]{0,0,0}\put(3981,-8912){\line( 5, 2){2381.897}}
}%
{\color[rgb]{0,0,0}\put(10860,-3431){\line( 5, 2){3295.862}}
}%
{\color[rgb]{0,0,0}\put(15510,-5634){\circle*{300}}
}%
{\color[rgb]{0,0,0}\put(14152,-5612){\circle*{300}}
}%
{\color[rgb]{0,0,0}\put(13177,-6579){\circle*{300}}
}%
{\color[rgb]{0,0,0}\put(13185,-7944){\circle*{300}}
}%
{\color[rgb]{0,0,0}\put(14130,-8897){\circle*{300}}
}%
{\color[rgb]{0,0,0}\put(11801,-4414){\circle*{266}}
}%
{\color[rgb]{0,0,0}\put(11496,-7936){\circle*{300}}
}%
{\color[rgb]{0,0,0}\put(11482,-6587){\circle*{300}}
}%
{\color[rgb]{0,0,0}\put(10530,-5619){\circle*{300}}
}%
{\color[rgb]{0,0,0}\put(9172,-5597){\circle*{300}}
}%
{\color[rgb]{0,0,0}\put(8197,-6564){\circle*{300}}
}%
{\color[rgb]{0,0,0}\put(8205,-7929){\circle*{300}}
}%
{\color[rgb]{0,0,0}\put(9150,-8882){\circle*{300}}
}%
{\color[rgb]{0,0,0}\put(10522,-8889){\circle*{300}}
}%
{\color[rgb]{0,0,0}\put(15502,-8904){\circle*{300}}
}%
{\color[rgb]{0,0,0}\put(13139,-4434){\circle*{266}}
}%
{\color[rgb]{0,0,0}\put(14122,-3451){\circle*{266}}
}%
{\color[rgb]{0,0,0}\put(14153,-2093){\circle*{266}}
}%
{\color[rgb]{0,0,0}\put(13180,-1140){\circle*{266}}
}%
{\color[rgb]{0,0,0}\put(11832,-1110){\circle*{266}}
}%
{\color[rgb]{0,0,0}\put(10859,-2073){\circle*{266}}
}%
{\color[rgb]{0,0,0}\put(10849,-3451){\circle*{266}}
}%
{\color[rgb]{0,0,0}\put(8502,-1341){\circle*{300}}
}%
{\color[rgb]{0,0,0}\put(7607,-2231){\circle*{300}}
}%
{\color[rgb]{0,0,0}\put(6405,-2236){\circle*{300}}
}%
{\color[rgb]{0,0,0}\put(5511,-1337){\circle*{300}}
}%
{\color[rgb]{0,0,0}\put(5506,-4323){\circle*{300}}
}%
{\color[rgb]{0,0,0}\put(6401,-3419){\circle*{300}}
}%
{\color[rgb]{0,0,0}\put(7598,-3415){\circle*{300}}
}%
{\color[rgb]{0,0,0}\put(8511,-4327){\circle*{300}}
}%
{\color[rgb]{0,0,0}\put(6290,-7944){\circle*{300}}
}%
{\color[rgb]{0,0,0}\put(6276,-6595){\circle*{300}}
}%
{\color[rgb]{0,0,0}\put(5324,-5627){\circle*{300}}
}%
{\color[rgb]{0,0,0}\put(3966,-5605){\circle*{300}}
}%
{\color[rgb]{0,0,0}\put(2991,-6572){\circle*{300}}
}%
{\color[rgb]{0,0,0}\put(2999,-7937){\circle*{300}}
}%
{\color[rgb]{0,0,0}\put(3944,-8890){\circle*{300}}
}%
{\color[rgb]{0,0,0}\put(5316,-8897){\circle*{300}}
}%
{\color[rgb]{0,0,0}\put(16476,-7951){\circle*{300}}
}%
{\color[rgb]{0,0,0}\put(16462,-6602){\circle*{300}}
}%
\put(9531,-9452){\makebox(0,0)[lb]{\smash{{\SetFigFont{20}{24.0}{\rmdefault}{\mddefault}{\updefault}{\color[rgb]{0,0,0}$H_2$}%
}}}}
\put(14721,-9437){\makebox(0,0)[lb]{\smash{{\SetFigFont{20}{24.0}{\rmdefault}{\mddefault}{\updefault}{\color[rgb]{0,0,0}$H_3$}%
}}}}
\put(6631,-4786){\makebox(0,0)[lb]{\smash{{\SetFigFont{20}{24.0}{\rmdefault}{\mddefault}{\updefault}{\color[rgb]{0,0,0}$C_4\times P_2$}%
}}}}
\put(4296,-9452){\makebox(0,0)[lb]{\smash{{\SetFigFont{20}{24.0}{\rmdefault}{\mddefault}{\updefault}{\color[rgb]{0,0,0}$H_1$}%
}}}}
\put(12309,-4912){\makebox(0,0)[lb]{\smash{{\SetFigFont{20}{24.0}{\rmdefault}{\mddefault}{\updefault}{\color[rgb]{0,0,0}$M_8$}%
}}}}
\end{picture}%
}
\caption{All cubic graphs of order $8$.}\label{fig-cubic8}
\end{center}
\end{figure}

\nt Clearly $\chi_{la}(K_{4})=4$. It was shown that $\chi_{la}(K_{3,3})=3$ in \cite{LNS} and $\chi_{la}(M_{2n})=3$ for odd $n$ in \cite{LNS-dmgt}. There is a conjecture proposed in \cite{LNS-dmgt} that $\chi_{la}(M_{2n})=4$ for even $n$.

\ms \nt We first consider $C_3\times P_2$. Note that $\chi_{la}(C_3\times P_2)\ge\chi(C_3\times P_2)=3$. Suppose there is a local antimagic $3$-coloring of $C_3\times P_2$ with colors $x,y,z$. By the uniqueness of 3-coloring of $C_3\times P_2$, we may let the labeling be shown as the figure below.

\centerline{\scalebox{0.4}{\begin{picture}(0,0)%
\includegraphics{C3xP2.pstex}%
\end{picture}%
\setlength{\unitlength}{3947sp}%
\begingroup\makeatletter\ifx\SetFigFont\undefined%
\gdef\SetFigFont#1#2#3#4#5{%
  \reset@font\fontsize{#1}{#2pt}%
  \fontfamily{#3}\fontseries{#4}\fontshape{#5}%
  \selectfont}%
\fi\endgroup%
\begin{picture}(5684,4990)(2637,-5342)
\put(5386,-661){\makebox(0,0)[lb]{\smash{{\SetFigFont{20}{24.0}{\rmdefault}{\mddefault}{\updefault}{\color[rgb]{0,0,0}$x$}%
}}}}
\put(8002,-5205){\makebox(0,0)[lb]{\smash{{\SetFigFont{20}{24.0}{\rmdefault}{\mddefault}{\updefault}{\color[rgb]{0,0,0}$y$}%
}}}}
\put(2768,-5186){\makebox(0,0)[lb]{\smash{{\SetFigFont{20}{24.0}{\rmdefault}{\mddefault}{\updefault}{\color[rgb]{0,0,0}$z$}%
}}}}
\put(5674,-3530){\makebox(0,0)[lb]{\smash{{\SetFigFont{20}{24.0}{\rmdefault}{\mddefault}{\updefault}{\color[rgb]{0,0,0}$e$}%
}}}}
\put(6237,-4112){\makebox(0,0)[lb]{\smash{{\SetFigFont{20}{24.0}{\rmdefault}{\mddefault}{\updefault}{\color[rgb]{0,0,0}$x$}%
}}}}
\put(5473,-4999){\makebox(0,0)[lb]{\smash{{\SetFigFont{20}{24.0}{\rmdefault}{\mddefault}{\updefault}{\color[rgb]{0,0,0}$c$}%
}}}}
\put(6927,-4787){\makebox(0,0)[lb]{\smash{{\SetFigFont{20}{24.0}{\rmdefault}{\mddefault}{\updefault}{\color[rgb]{0,0,0}$h$}%
}}}}
\put(3790,-4387){\makebox(0,0)[lb]{\smash{{\SetFigFont{20}{24.0}{\rmdefault}{\mddefault}{\updefault}{\color[rgb]{0,0,0}$i$}%
}}}}
\put(5566,-1817){\makebox(0,0)[lb]{\smash{{\SetFigFont{20}{24.0}{\rmdefault}{\mddefault}{\updefault}{\color[rgb]{0,0,0}$g$}%
}}}}
\put(5379,-2716){\makebox(0,0)[lb]{\smash{{\SetFigFont{20}{24.0}{\rmdefault}{\mddefault}{\updefault}{\color[rgb]{0,0,0}$z$}%
}}}}
\put(5323,-3931){\makebox(0,0)[lb]{\smash{{\SetFigFont{20}{24.0}{\rmdefault}{\mddefault}{\updefault}{\color[rgb]{0,0,0}$f$}%
}}}}
\put(5078,-3530){\makebox(0,0)[lb]{\smash{{\SetFigFont{20}{24.0}{\rmdefault}{\mddefault}{\updefault}{\color[rgb]{0,0,0}$d$}%
}}}}
\put(4184,-3076){\makebox(0,0)[lb]{\smash{{\SetFigFont{20}{24.0}{\rmdefault}{\mddefault}{\updefault}{\color[rgb]{0,0,0}$a$}%
}}}}
\put(6587,-3105){\makebox(0,0)[lb]{\smash{{\SetFigFont{20}{24.0}{\rmdefault}{\mddefault}{\updefault}{\color[rgb]{0,0,0}$b$}%
}}}}
\put(4607,-4101){\makebox(0,0)[lb]{\smash{{\SetFigFont{20}{24.0}{\rmdefault}{\mddefault}{\updefault}{\color[rgb]{0,0,0}$y$}%
}}}}
\end{picture}%
}}

\nt It is clear that \begin{align}2(a+b+c+d+e+f+g+h+i)=& 2(x+y+z)=2\times \frac{9\times 10}{2}=90\nonumber\\
2(a+b+c)+(g+h+i)=& x+y+z=45\label{eq-2}\\
2(d+e+f)+(g+h+i)= & x+y+z\nonumber\end{align}
So we have $a+b+c=d+e+f$. Since $a+c+i=z=e+d+g$ and $e+f+h=x=a+b+g$, $a+i=e+d+g-c$ and $e+h=a+b+g-f$. Since $a+b-f=d+e-c$, $a+i=e+h$. Similarly, we have $b+g=f+i$ and $c+h=d+g$.

\nt From \eqref{eq-2} we have $21\le 2(a+b+c)\le 39$, i.e., $11\le a+b+c\le 19$. Without loss of generality, we may assume that $\min\{a,b,c,d,e,f\}=f$.
\begin{enumerate}[1.]
\item Suppose $a+b+c=d+e+f=19$. There are 5 combinations of three distinct labels whose sum is 19, which are $\{9,8,2\}$, $\{9,7,3\}$, $\{9,6,4\}$, $\{8,7,4\}$, $\{8,6,5\}$. Since $a,b,c,d,e,f$ are distinct, only one pair of combinations can be chosen which is $\{\{a,b,c\}, \{d,e,f\}\}=\{\{9,7,3\}, \{8,6,5\}\}$. By the assumption, $f=3$, $b\in \{8,6,5\}$ and $\{g,h,i\}=\{1,2,4\}$. Then $6\le b+g=f+i\in\{4,5,7\}$. Thus $i=4$ and $(b,g)=(6,1), (5,2)$. Two possible solutions are $(a,b,c,d,e,f,g,h,i)=(5,6,8,9,7,3,1,2,4)$ and $(6,5,8,7,9,3,2,1,4)$.

\item Suppose $a+b+c=d+e+f=18$. There are 7 combinations: $\{9,8,1\}$, $\{9,7,2\}$, $\{9,6,3\}$, $\{9,5,4\}$, $\{8,7,4\}$, $\{8,6,4\}$, $\{7,6,5\}$. Only two cases of $\{\{a,b,c\}, \{d,e,f\}\}$ need to be considered. Under the assumption we have:
    \begin{enumerate}[2.1)]
    \item $\{a,b,c\}=\{7,6,5\}$ and $\{d,e,f\}=\{9,8,1\}$, where $f=1$. Hence $\{g,h,i\}=\{2,3,4\}$. Now $b+g=f+i\in\{3,4,5\}$. There is no solution.
    \item $\{a,b,c\}=\{8,6,4\}$ and $\{d,e,f\}=\{9,7,2\}$, where $f=2$. Hence $\{g,h,i\}=\{1,3,5\}$. Note that $4+g\le b+g=f+i\in\{3,5,7\}$. This implies that $i\ge 3$ and $g\le 3$.

         If $i=3$, then $g=1$. Hence $b=4$, $h=5$. But $e+h=a+i\in\{11, 9\}$ is impossible.

        If $i=5$, then $b=6$ and $g=1$, or $b=4$ and $g=3$.
         \begin{enumerate}[{2.2}.1)]
        \item If $b=6$ and $g=1$, then $h=3$. But $e+h=a+i\in\{13, 9\}$ is impossible.
        \item If $b=4$ and $g=3$, then $h=1$. But $e+h=a+i\in\{13, 11\}$ is impossible.
        \end{enumerate}
    \end{enumerate}
    \end{enumerate}
\nt For the cases of $a+b+c=d+e+f\in[11,17]$, the argument are similar. We omit them. We only list all the results as follows:\\
    \renewcommand{\multirowsetup}{\centering}
    $\begin{array}{c|c|l}
    \mbox{Case} & a+b+c & (a,b,c,d,e,f,g,h,i)\\\hhline{===}
    1 & \multirow{2}{1cm} {19} & (5,6,8,9,7,3,1,2,4)\\
    2 & & (6,5,8,7,9,3,2,1,4)\\\hline
    3 & \multirow{5}{1cm} {17} & (8,3,6,9,7,1,2,5,4)\\
    4 & & (5,4,8,7,9,1,3,2,6)\\
    5 & & (4,5,8,9,7,1,2,3,6)\\
    6 & & (5,8,4,6,9,2,1,3,7)\\
    7 & & (9,5,3,8,7,2,1,6,4)\\\hline
    8 & \multirow{3}{1cm} {16} & (4,7,5,6,9,1,2,3,8)\\
    9 & & (9,2,5,8,7,1,3,6,4)\\
    10 & & (9,4,3,7,8,1,2,6,5)\\\hline
    \end{array}$ \qquad
    $\begin{array}{c|c|l}
    \mbox{Case} & a+b+c & (a,b,c,d,e,f,g,h,i)\\\hhline{===}
    11 & \multirow{3}{1cm} {14} & (5,6,3,9,4,1,2,8,7)\\
    12 & & (9,3,2,5,8,1,4,7,6)\\
    13 & & (9,2,3,7,6,1,4,8,5)\\\hline
    14 & \multirow{5}{1cm} {13}& (7,4,2,3,9,1,5,6,8)\\
    15 & & (2,5,6,9,3,1,4,7,8)\\
    16 & & (5,2,6,3,9,1,7,4,8)\\
    17 & & (6,5,2,8,4,1,3,9,7)\\
    18 & & (8,3,2,7,5,1,4,9,6)\\\hline
    \end{array}$\\

\nt Thus, up to symmetry, there are 18 local antimagic $3$-colorings of $C_3\times P_2$. Following only shows the result obtained at Case~6.

\centerline{\epsfig{file=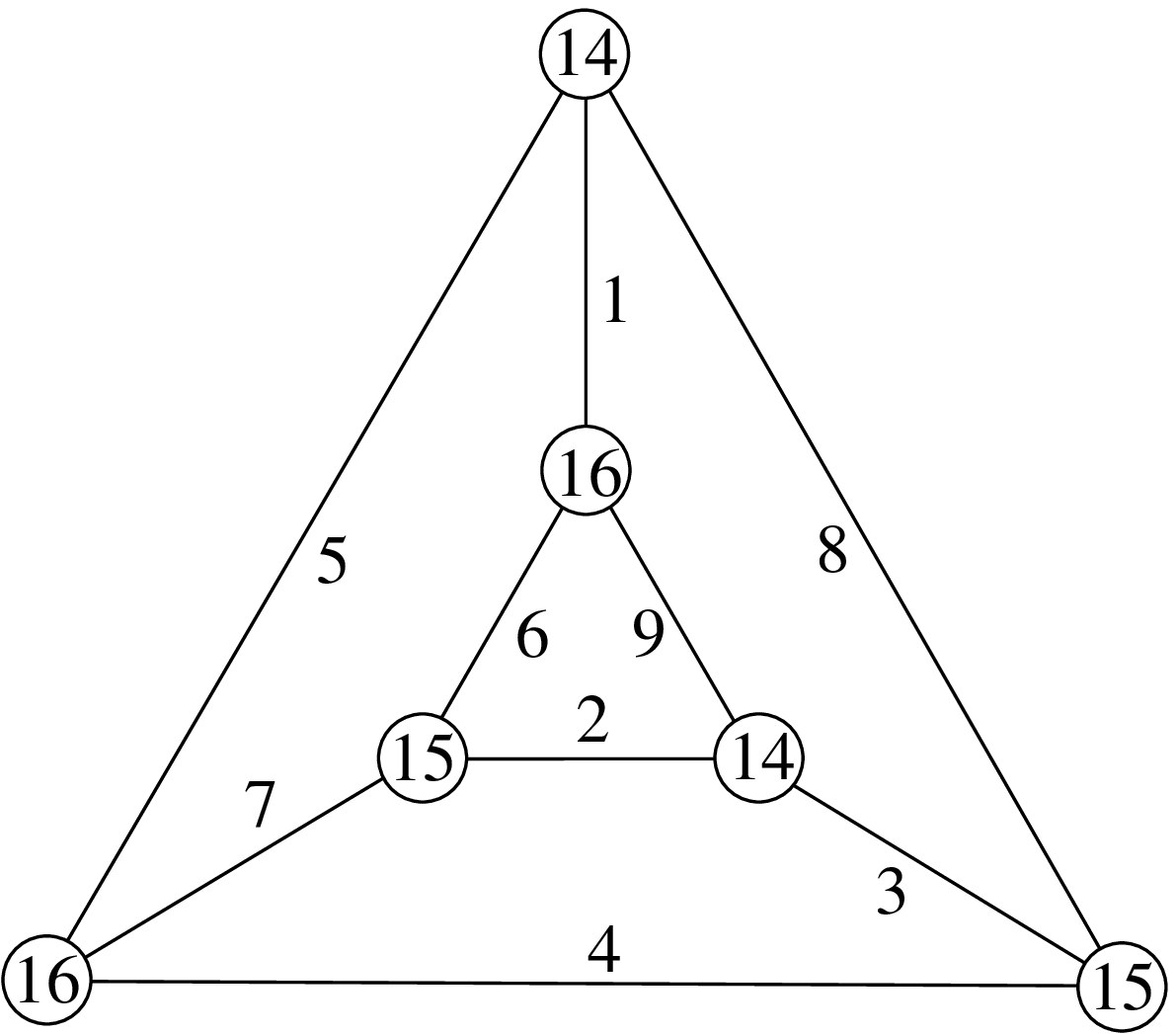, width=4cm}}

\nt From the above discussion, and by Lemma~\ref{lem-reg}, we have
\begin{theorem}
$\chi_{la}(C_3\times P_2) = \chi_{la}((C_3\times P_2)-e) =3$. Moreover, up to symmetry, $C_3\times P_2$ admits $18$ different local antimagic $3$-colorings.
\end{theorem}

\begin{lemma}\label{lem-4-4bipartite}
Let $G=(X,Y)$ be a bipartite graph of size $q$ with $|X|=|Y|=4$. If $C_4\times P_2$ is a spanning subgraph of $G$ and $12\le q\le 14$, then $\chi_{la}(G)\ge 4$.\end{lemma}

\begin{proof} Since $|X|=|Y|$,  $\chi_{la}(G)\ge 3$ by Corollary~\ref{cor-2part}. Suppose there were a local antimagic $3$-coloring of $G$. Let $V_a, V_b, V_c$ be the induced color classes such that each vertex in $V_j$ is colored by $j$, $j\in\{a,b,c\}$.

\ms\nt Suppose a color set containing at least 4 vertices, say $V_a$. Since each vertex in $X$ is adjacent to at least 3 vertices in $Y$ and vice versa, $X\subseteq V_a$ or $Y\subseteq V_a$.
Then $\sum\limits_{x\in X}f^+(x)=4a=\frac{q(q+1)}{2}$. Thus $q(q+1)\equiv 0\pmod 8$. But this equation does not hold when $12\le q\le 14$.  Thus, $|V_j|\le 3$.

\ms\nt Without loss of generality, we assume $|V_a|=3=|V_b|$ and $|V_c|=2$. Since each vertex in $X$ is adjacent to at least 3 vertices in $Y$ and vice versa, we may assume that $V_a\subset X$ and $V_b\subset Y$; also one vertex in $V_c$ lying in $X$ and the other one lying in $Y$. Thus $3a+c=3b+c=\frac{q(q+1)}{2}$. Hence $a=b$ which is a contradiction.
\end{proof}

\begin{theorem}\label{thm-C4P2}
$\chi_{la}(C_4\times P_2)=\chi_{la}((C_4\times P_2)+e)=4$, where $e$ is an extra edge whose end vertices lie in different partition of $C_4\times P_2$.
\end{theorem}
\begin{proof} By Lemma~\ref{lem-4-4bipartite}, both $\chi_{la}(C_4\times P_2)$ and $\chi_{la}(C_4\times P_2+e)$ are at least $4$.

\nt Actually there is a local antimagic $4$-coloring of $C_4\times P_2$ which is shown below:\\[1mm]
\centerline{\epsfig{file=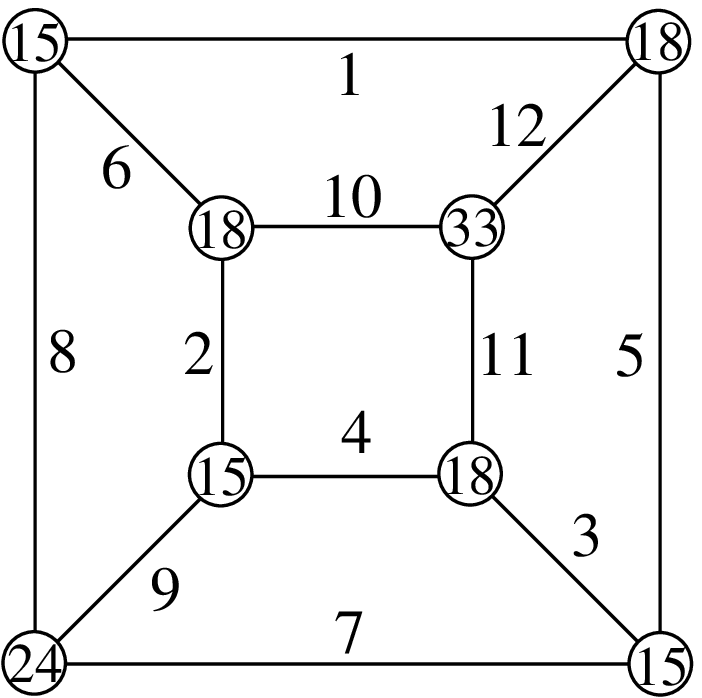, width=3cm}}

\nt Thus $\chi_{la}(C_4\times P_2)=4$.

\ms\nt It is easy to see that $(C_4\times P_2)+e$ is unique up to isomorphism if the end vertices of $e$ lie in different partition. We have a local antimagic $4$-coloring of $(C_4\times P_2)+e$ which is the labeling obtained from the above labeling of $C_4\times P_2$ by adding an extra edge with label 13 joining the vertices with label 24 and 33. Thus $\chi_{la}((C_4\times P_2)+e)=4$.
\end{proof}

\nt Now let us consider a cubic graph $G$ of order 8.

\begin{theorem}\label{thm-cubic8} Suppose $G$ is a cubic graph of order $8$, then $\chi_{la}(G)=4$.
\end{theorem}
\begin{proof} By Theorem~\ref{thm-C4P2}, we only need to consider $G\in\{M_8, H_1, H_2, H_3\}$. We follow the notation in Figure~\ref{fig-cubic8}.
Clearly, $\chi_{la}(G)\ge \chi(G)=3$. Suppose there were a local antimagic $3$-coloring $f$ of $G$. Let $V_x, V_y, V_z$ be the corresponding color classes so that $f^+(u)=j$ if $u\in V_j$. Since the independent number of $G$ is 3, $|V_j|\le 3$. Without loss of generality, we let $V_x=\{u_1, u_2\}$, $V_y=\{v_1,v_2,v_3\}$ and $V_z=\{w_1, w_2, w_3\}$. Let $E_j$ be the set of edges incident with vertices in $V_j$, $j\in\{x,y,z\}$. Clearly,$|E_x|=6$ and  $|E_y|=|E_z|=9$.

\ms\nt We may assume that edges with labels $a_1, a_2, a_3$ are incident to $v_1$, edges with labels $a_4, a_5, a_6$ are incident to $v_2$ and edges with labels $a_7, a_8, a_9$ are incident to $v_3$. Note that, another 3 edges not incident to vertices in $V_y$ must be incident with vertices in $V_x$.

\ms\nt Similarly, we assume that edges with labels $b_1, b_2, b_3$ are incident to $u_1$, edges with labels $b_4, b_5, b_6$ are incident to $u_2$. Thus we have
\begin{align} x& =  b_1+b_2+b_3=b_4+b_5+b_6\label{eq-x}
\\
y & = a_1+a_2+a_3 =a_4+a_5+a_6=a_7+a_8+a_9\label{eq-y}\\
2x+3y+3z & = 2\times\mbox{(sum of all labels)}=132\nonumber
\end{align}
We will also have a formula about $z$ similar to $y$.

\ms\nt From \eqref{eq-x} we have $21=\sum\limits_{i=1}^{6} i\le 2x=\sum\limits_{i=1}^6 b_i\le \sum\limits_{i=7}^{12} i=57$. So
$22\le 2x\le 56$. Similarly we have $45\le 3y, 3z\le 72$. It is equivalent to
\[x \in[11,28] \mbox{ and } y,z \in[15,24].\]
Since $22\le 2x\le 56$, $76=132-56\le 3y+3z\le 132-22=110$. Hence $26\le y+z\le 36$.
Since $2x+3y+3z=132$ and $x,y,z$ are distinct, we have the following 5 cases for $(x, \{y,z\})$:
 \[(18, \{15, 17\}), \ (15, \{16, 18\}),\ (12, \{15, 21\}), \ (12, \{16, 20\}),\ (12, \{17, 19\})
\]

\ms\nt According to \eqref{eq-y}, we have to find a set of 3 (disjoint) triples, namely\\ $Y=\{\{a_1,a_2,a_3\},\{a_4,a_5,a_6\}, \{a_7,a_8, a_9\}\}$, such that the sum of labels in every triple is $y$. Similarly, we have to find a set of 3 (disjoint) triples such that the sum of labels in every triple is $z$.

\ms \nt Similarly, according to \eqref{eq-x}, we have to find a set of 2 (disjoint) triples, namely\\ $X=\{\{b_1,b_2,b_3\}, \{b_4,b_5,b_6\}\}$, such that the sum of labels in every triple is $x$. Moreover, 3 labels not appear in $Y$ (or $Z$) must lie in $X$.

\ms \nt Let us consider the combinations of a set of 3 triples with distinct labels satisfying \eqref{eq-y} so that the sum of labels in each triple is $s$, where $15\le s\le 17$:
\begin{enumerate}[1.]
\item Suppose $s=15$. There are 12 combinations. They are $\{12,2,1\}$, $\{11,3,1\}$, $\{10,4,1\}$, $\{10,3,2\}$, $\{9,5,1\}$, $\{9,4,2\}$, $\{8,6,1\}$, $\{8,5,2\}$, $\{8,4,3\}$, $\{7,6,2\}$, $\{7,5,3\}$, $\{6,5,4\}$. There are two possible sets of 3 triples satisfying \eqref{eq-y}: $\{\{9,5,1\},\{8,4,3\}, \{7,6,2\}\}$ and $\{\{9,4,2\}, \{8,6,1\}, \{7,5,3\}\}$.

\item Suppose $s=16$. Similar to Case~1, there are 13 combinations: $\{12,3,1\}$, $\{11,4,1\}$, $\{11,3,2\}$, $\{10,5,1\}$, $\{10,4,2\}$, $\{9,6,1\}$, $\{9,5,2\}$, $\{9,4,3\}$, $\{8,7,1\}$, $\{8,6,2\}$, $\{8,5,3\}$, $\{7,6,3\}$, $\{7,5,4\}$. There are four solutions:\\
$\{\{11,4,1\}, \{9,5,2\}, \{7,6,3\}\}$, $\{\{11,3,2\}, \{9,6,1\}, \{7,5.4\}\}$, $\{\{10,5,1\}, \{9,4,3\}, \{8,6,2\}\}$ and\break $\{\{10,4,2\}, \{9,6,1\}, \{8,5,3\}\}$.

\item Suppose $s=17$. There are 14 combinations: $\{12,4,1\}$, $\{12,3,2\}$, $\{11,5,1\}$, $\{11,4,2\}$, $\{10,6,1\}$, $\{10,5,2\}$, $\{10,4,3\}$,  $\{9,7,1\}$, $\{9,6,2\}$, $\{9,5,3\}$, $\{8,7,2\}$, $\{8,6,3\}$, $\{8,5,4\}$, $\{7,6,4\}$. There are seven solutions:\\
$\{\{12,4,1\},\{10,5,2\}, \{8,6,3\}\}$,  $\{\{12,3,2\},\{11,5,1\}, \{7,6,4\}\}$, $\{\{11,5,1\}, \{10,4,3\}, \{9,6,2\}\}$,\\ $\{\{11,5,1\}, \{10,4,3\}, \{8,7,2\}\}$, $\{\{11,4,2\}, \{10,6,1\}, \{9,5,3\}\}$, $\{\{11,4,2\}, \{9,7,1\}, \{8,6,3\}\}$ and\break  $\{\{10,5,2\}, \{9,7,1\}, \{8,6,3\}\}$

\end{enumerate}

\begin{enumerate}[a.]
\item Suppose $y=15$ (we may apply the same argument for $z=15$). From Case~1, there are two possible sets of 3 triples for $Y$ which are $\{\{9,5,1\},\{8,4,3\}, \{7,6,2\}\}$ and $\{\{9,4,2\}, \{8,6,1\}, \{7,5,3\}\}$. Since 12, 11 and 10 do not appear in $Y$, they appear in $X$. At least two of them lie in a triple of $X$. Then $x> 21$, which is not a case.

\item Suppose $y=16$ (we may apply the same argument for $z=16$). From Case~2, there are four cases:\\
$\{\{11,4,1\}, \{9,5,2\}, \{7,6,3\}\}$, $\{\{11,3,2\}, \{9,6,1\}, \{7,5.4\}\}$, $\{\{10,5,1\}, \{9,4,3\}, \{8,6,2\}\}$ and\break $\{\{10,4,2\}, \{9,6,1\}, \{8,5,3\}\}$. Here, either 12,10,8 or 12,11,7 do not appear in $Y$. Similar to Case~a, we have $x>18$ which is not a case.

\item Suppose $y=17$ (we may apply the same argument for $z=17$). From Case~3 and by a similar argument as the above cases, we may get $x>15$. But the remaining case is $(x, \{y,z\})=(12, \{17, 19\})$.
\end{enumerate}
\nt Thus, there is no local antimagic $3$-coloring of $G$. Hence $\chi_{la}(G)\ge 4$. A local antimagic 4-coloring of each $G$ is shown below:\\[1mm]
\centerline{\epsfig{file=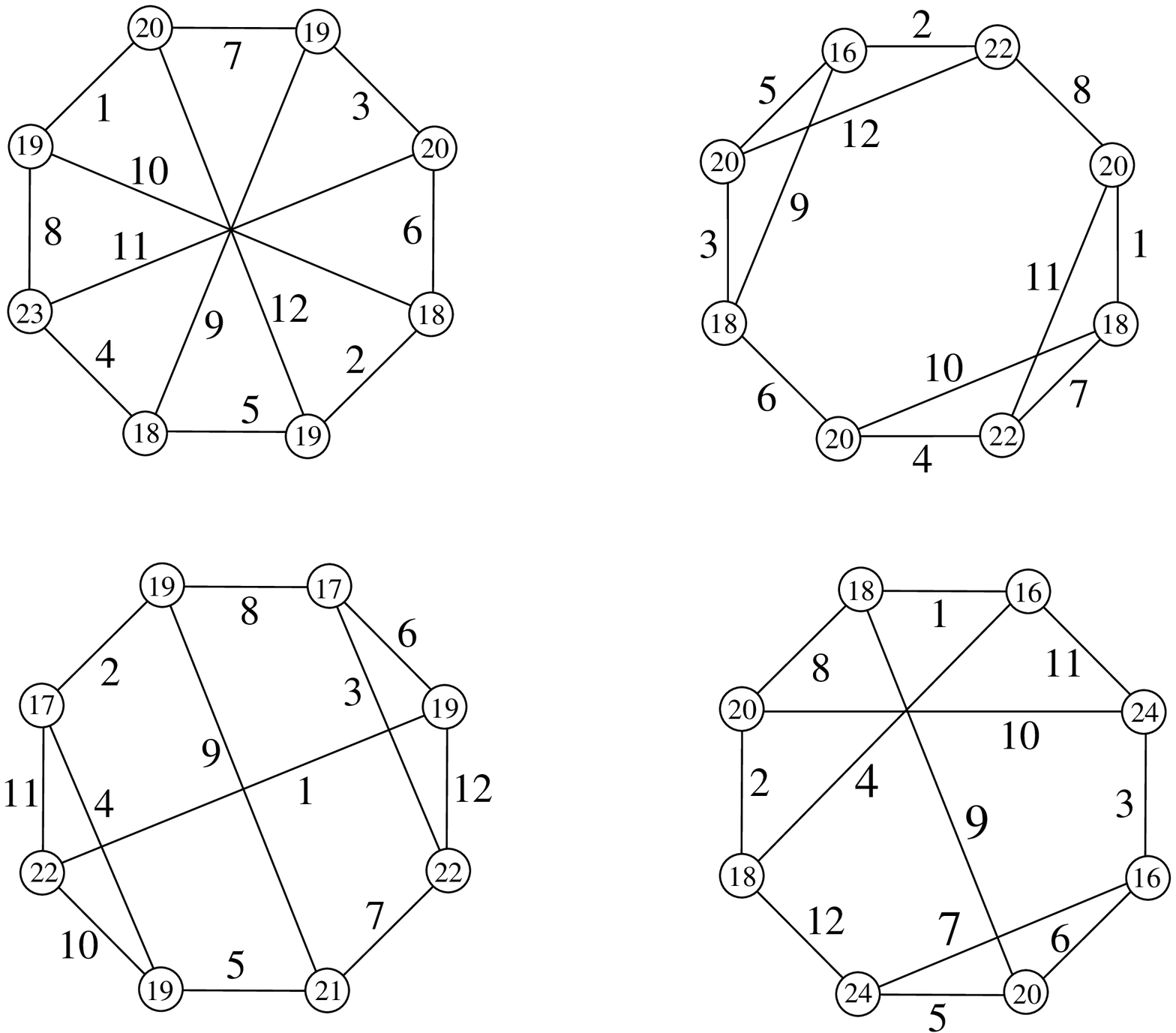, width=9cm}}

\nt Thus $\chi_{la}(G)=4$. \end{proof}

\subsection{Connected Quartic Graphs of Order at Most 8}\label{sec-quartic}

In this subsection, we consider connected quartic graphs or their related graphs of order at most 8. There is only one quartic graph of order 5 which is $K_5$. It is known that $\chi_{la}(K_5)=5$. There is only one quartic graph of order 6 which is the octahedral graph, say $Q$. Clearly, $\chi_{la}(Q)\ge\chi(Q)=3$. The local antimagic labeling of $G$ given below implies that $\chi_{la}(Q)=3$.

\centerline{\epsfig{file=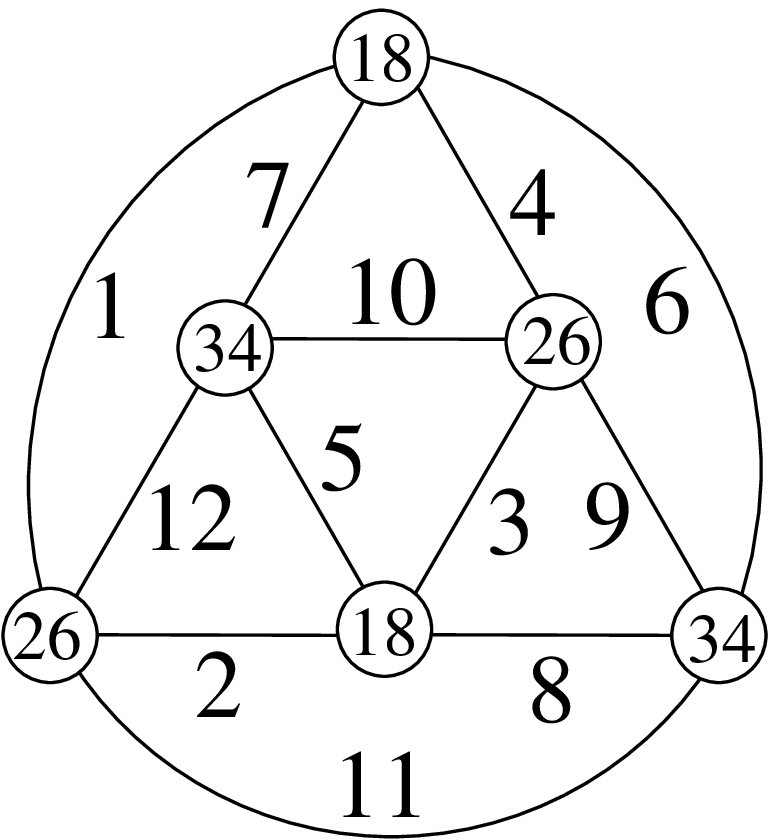, width=2.5cm}}

\nt There are two quartic graphs of order 7, which are $C_7(1,2)$ and $\overline{C_3+C_4}$, the complement of $C_3+C_4$. A local antimagic labeling for each graph is as shown below.

\vspace{0.5cm}\centerline{\epsfig{file=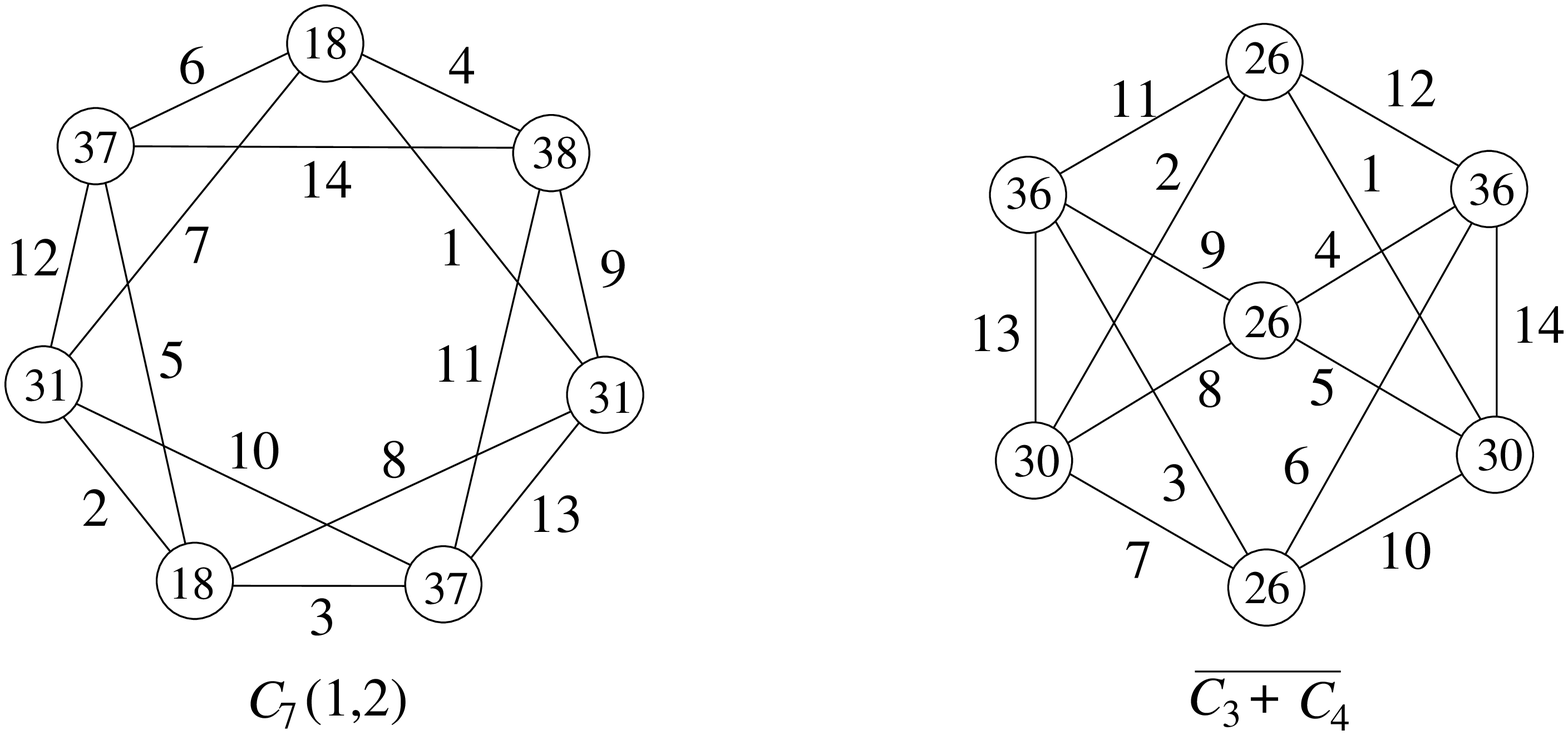, width=8cm}}

\ms\nt It is easy to see that $\chi(C_7(1,2))=4$ and $\chi(\overline{C_3+C_4})=3$. So $\chi_{la}(C_7(1,2))=4$ and $\chi_{la}(\overline{C_3+C_4})=3$.

\ms \nt It is known that there are 6 quartic graphs of order 8. They are $B_1=K_{4,4}$ which is known that $\chi_{la}(K_{4,4})=3$, $B_2=K_4\times K_2$, $B_3=C_8(1,2)$, and the other three are denoted by $B_4$, $B_5$, $B_6$ as shown. Local antimagic labelings for $B_i, 2\le i\le 6$ are also given.

\vspace{0.5cm}\centerline{\epsfig{file=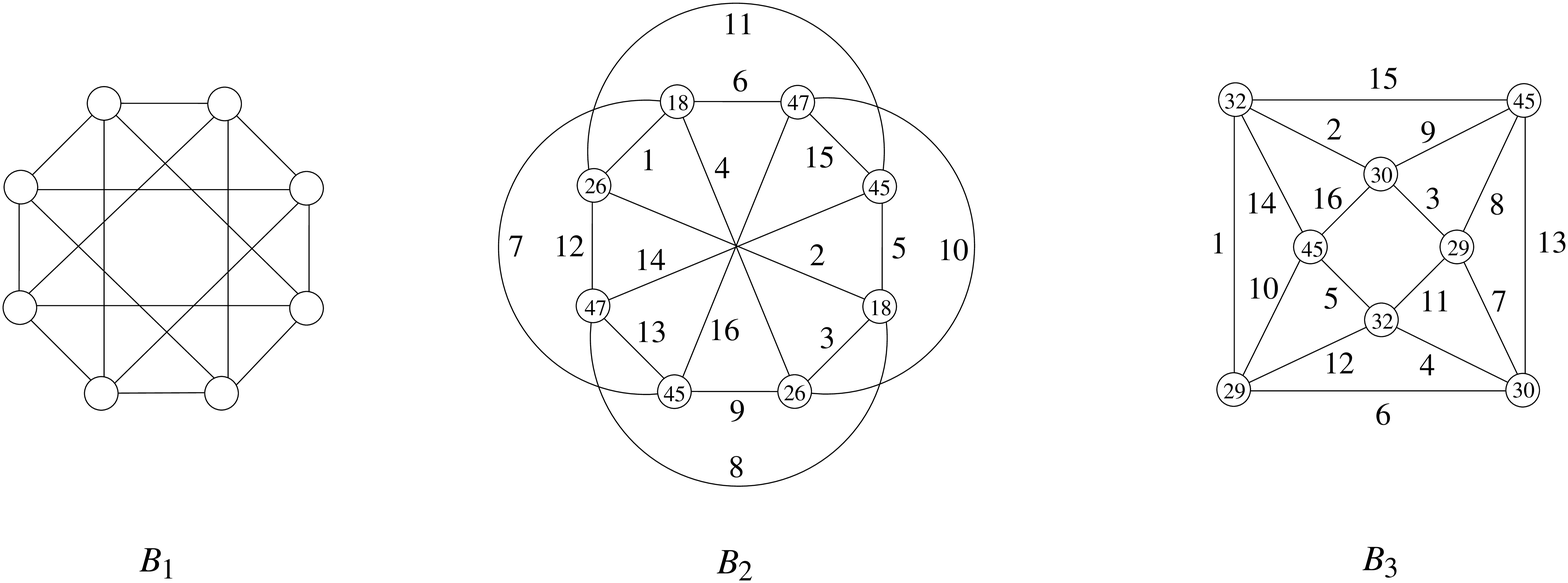, width=16cm}}

\medskip

\vspace{0.5cm}\centerline{\epsfig{file=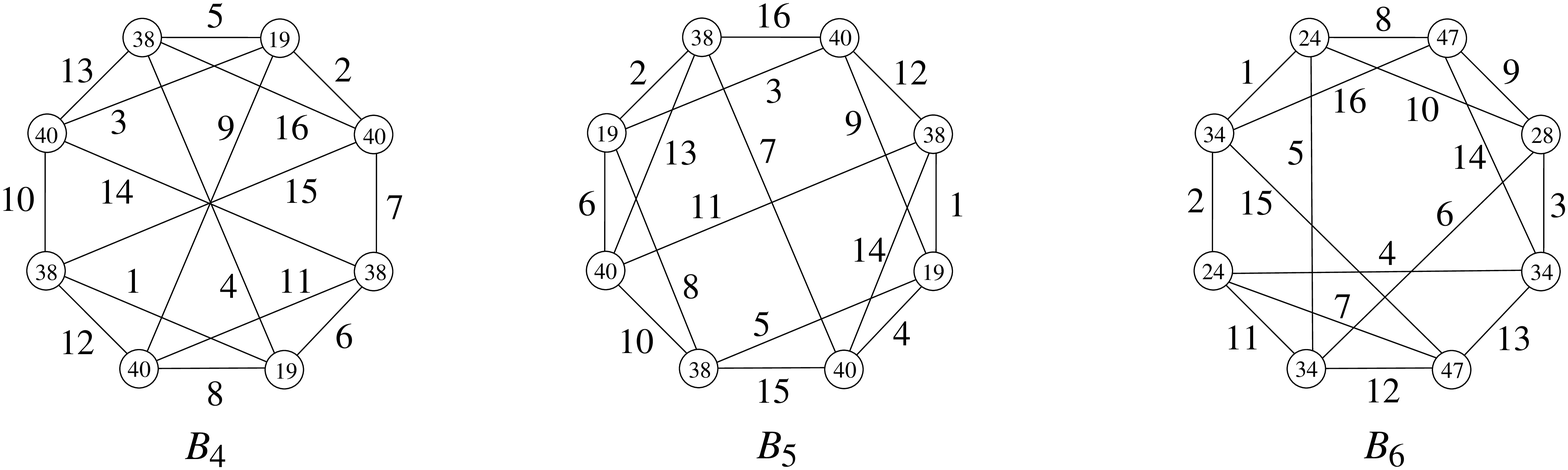, width=15cm}}

\ms\nt Clealry, $\chi(B_2)=4$. It can be showed that $\chi(B_4)=\chi(B_5)=3$, $\chi(B_3)=\chi(B_6)=4$.
Thus $\chi_{la}(B_2)=4$, $\chi_{la}(B_3)=4$, $\chi_{la}(B_4)=3$, $\chi_{la}(B_5)=3$ and $\chi_{la}(B_6)=4$. Combining the above discussion we have

\begin{theorem} Suppose $G$ is a quartic of order $p$, where $5\le p\le 8$, then $\chi_{la}(G)=\chi(G)$ except $\chi_{la}(K_{4,4})=\chi(K_{4,4})+1$.
\end{theorem}

\subsection{Connected Quintic Graphs of Order at Most 8}\label{sec-quintic}

\nt  There is one quintic graph of order 6 which is $K_6$. Obviously, $\chi_{la}(K_6)=6$. There are 3 quintic graphs of order 8: $C_8(1,2,4)$, $C_8(1,3,4)\cong 2K_2 \vee 2K_2$ and $(5,3)$-cone, as shown below.

\vskip0.5cm\centerline{\epsfig{file=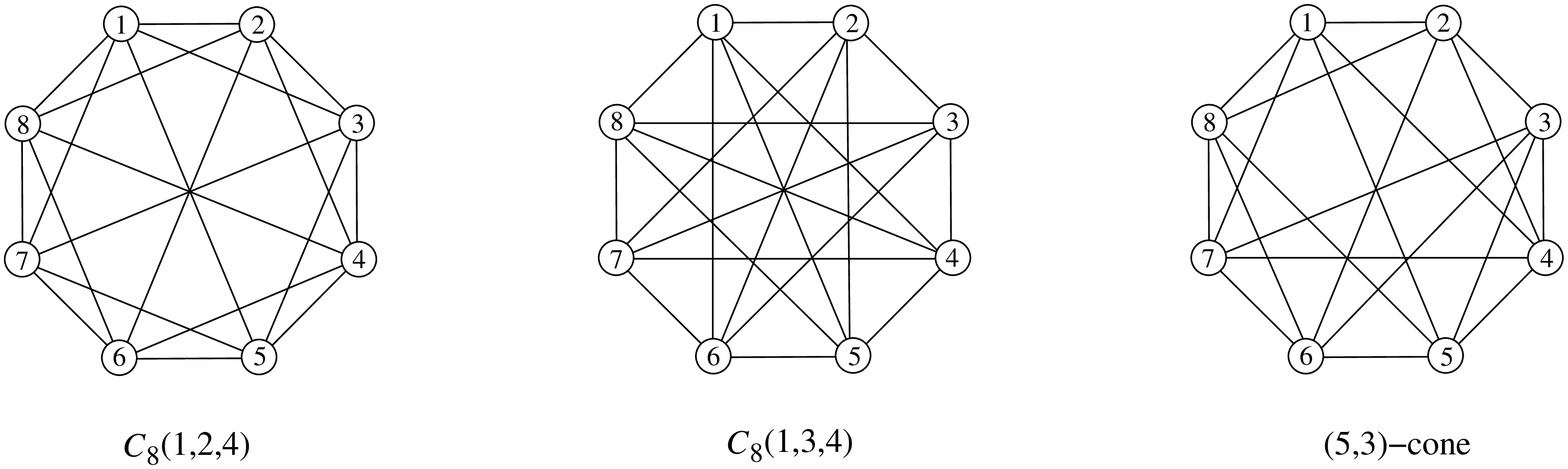, width=12cm}}

\nt We use $\{1,2,3,4,5,6,7,8\}$ as vertex names showed in the above figure. Note that $(5,3)$-cone is isomorphic to $C_5\vee O_3$. Namely, the $5$-cycle is $(14368)$. It was showed in \cite{LNS-dmgt} that $\chi_{la}(C_5\vee O_3)=4$. Here are labeling matrices of local antimagic labelings for $C_8(1,2,4)$ and $C_8(1,3,4)$, respectively. The last column shows its row sums.
\begin{align*}M(C_8(1,2,4))&=\left(\begin{array}{*{8}{c}}
* & 19 & 12 & * & 1 & * & 20 & 11\\
19 & * & 6 & 17 & * & 5 & * & 9\\
12 & 6 & * & 2 & 18 & * & 4 & *\\
* & 17 & 2 & * & 15 & 13 & * & 16\\
1 & * & 18 & 15 & * & 7 & 8 & *\\
* & 5 & * & 13 & 7 & * & 14 & 3\\
20 & * & 4 & * & 8 & 14 & * & 10\\
11 & 9 & * & 16 & * & 3 & 10 & *
\end{array}\right)\begin{array}{c}63\\56\\42\\63\\49\\42\\56\\49\end{array},\\
M(C_8(1,3,4))&=\left(\begin{array}{*{8}{c}}
* & 17 & * & 1 & 20 & 3 & * & 7\\
17 & * & 9 & * & 10 & 13 & 2 & *\\
* & 9 & * & 12 & * & 16 & 19 & 11\\
1 & * & 12 & * & 14 & * & 18 & 6\\
20 & 10 & * & 14 & * & 8 & * & 15\\
3 & 13 & 16 & * & 8 & * & 4 & *\\
* & 2 & 19 & 18 & * & 4 & * & 5\\
7 & * & 11 & 6 & 15 & * & 5 & *
\end{array}\right)\begin{array}{c}48\\51\\67\\51\\67\\44\\48\\44\end{array}. \end{align*}
\nt Thus we have the following theorem.

\begin{theorem} Suppose $G$ is a quintic graph of order $8$, then $\chi_{la}(G)=\chi(G)=4$.
\end{theorem}

\subsection{Connected Sextic and Septic Graphs of Order at Most 8}\label{sec-sextic}

There is only one sextic graph of order 7 which is $K_7$. There is only one sextic graph of order 8 which is $C_8(1,2,3)$. There is only one septic graph of order 8 which is $K_8$.
\begin{theorem} $\chi_{la}(C_8(1,2,3))=\chi(C_8(1,2,3))=4$.
\end{theorem}
\begin{proof} Following is the labeling matrix of a local antimagic 4-coloring of $C_8(1,2,3)$.

\centerline{
$\left(\begin{array}{*{8}{c}}
* & 14 & 22 & 13 & * & 2 & 3 & 23\\
14 & * & 5 & 10 & 12 & * & 1 & 15\\
22 & 5 & * & 17 & 6 & 4 & * & 20\\
13 & 10 & 17 & * & 19 & 9 & 24 & *\\
* & 12 & 6 & 19 & * & 8 & 21 & 11\\
2 & * & 4 & 9 & 8 & * & 18 & 16\\
3 & 1 & * & 24 & 21 & 18 & * & 7\\
23 & 15 & 20 & * & 11 & 16 & 7 & *
\end{array}\right)\begin{array}{c}77\\57\\74\\92\\77\\57\\74\\92\end{array}$.
}
\end{proof}
\nt Let us summarize the local antimagic chromatic numbers of connected $r$-regular graphs of order at most 8:
\renewcommand{\arraystretch}{1.2}
\renewcommand{\multirowsetup}{\centering}
\[\begin{array}{|c||*{6}{c|}}\hline
\multirow{2}{1cm}{$r$} & \multicolumn{6}{c|}{\mbox{order}}\\\cline{2-7}
 & 3 & 4 & 5 & 6 & 7 & 8\\\hhline{|=::======|}
2 & \multicolumn{6}{c|}{3}\\\cline{1-7}
3 & * & 4 & * & 3 & * & 4\\\cline{1-7}
4 & * & * & 5 & 3 & \begin{matrix}\chi_{la}(C_7(1,2))=4,\\\chi_{la}(\overline{C_3+C_4})=3\end{matrix} & \begin{matrix}\chi_{la}(B_i)=4, i=2,3,6\\ \chi_{la}(B_j)=3, j=1,4,5\end{matrix}\\\cline{1-7}
5 & * & * & * & 6 & * & 4 \\\cline{1-7}
6 & * & * & * & * & 7 & 4\\\cline{1-7}
7 & * & * & * & * & * & 8\\\cline{1-7}
\end{array}\]

\section{Join Product of Regular Graphs}\label{sec-join}

\nt Let $G$ and $H$ be $r_1$-regular and $r_2$-regular with order $p_1$ and $p_2$, respectively. Then $G\vee H$ is an $r$-regular graph if $r_1 + p_2 = r_2 + p_1 = r$. The following two theorems on $G\vee O_n$ were obtained in~\cite{LNS}.

\begin{theorem}\label{thm-GVOn} Suppose $G$ is of order $m\ge 3$ with $m \equiv n\pmod{2}$ and $\chi(G)=\chi_{la}(G)$. If
(i) $n\ge m$, or (ii) $m\ge n^2/2$ and $n\ge 4$, then $\chi_{la}(G\vee O_n)=\chi_{la}(G)+1$. \end{theorem}

\begin{theorem}\label{thm-regGVOn} Suppose $G$ is an $r$-regular graph of order $m\ge 3$ with $m \equiv n\pmod{2}$ and $\chi(G)=\chi_{la}(G)$. If $m >n\ge 2$ and $r \ge \frac{(m-n)(mn+1)}{2mn}$, then $\chi_{la}(G\vee O_n)=\chi_{la}(G)+1$. \end{theorem}

\nt Suppose $G$ is a $2r$-regular graph of order $m$. If $n=m-2t$ $(t\ge 1)$, then $m\equiv n \pmod{2}$. Moreover, $G\vee O_{m-2r}$ is an $m$-regular graph of order $2m-2r$.  Combining Theorem~\ref{thm-regGVOn}, we immediately have the following theorem.

\begin{theorem}\label{thm-join-regular-On} Suppose $G$ is a $2r$-regular graph of order $m\ge 3$ and $\chi(G)=\chi_{la}(G)$. If $r\ge 2$, $r\ge t\ge 1$ and $m-2t\ge 2$, then $\chi_{la}(G\vee O_{m-2t}) = \chi_{la}(G) + 1$. Particularly, let $G_0=G$ and $G_{k+1} = G_{k} \vee O_{m-2r}$ for $k\ge 0$. If $m$ is even and $m-2r\ge 2$, then $G_{k+1}$ is a $(2r+k(m-2r))$-regular graph of order $m+(k+1)(m-2r)$ with $\chi_{la}(G_k)=\chi_{la}(G)+k$.    \end{theorem}

\nt In this section, we shall make use of some magic rectangles (for the existence, see~\cite{Hagedorn}).

\begin{lemma}\label{lem-col-row}
Suppose $M$ is an $S\times T$ magic rectangle using integers in $[Q+1, Q+ST]$. Let $C$ and $R$ be the column sum and row sum of $M$, respectively. Then
\[C-R =\frac{1}{2}(2Q+1+ST)(S-T).\]
\end{lemma}
\begin{proof}
$\displaystyle C-R =\frac{1}{2}(2Q+1+ST)(ST)\left[\frac{1}{T}-\frac{1}{S}\right] =\frac{1}{2}[2Q+1+ST](S-T).$
\end{proof}

\nt Suppose $A=(a_1,\dots, a_m)$ and $B=(b_1,\dots, b_n)$ are sequences of length $m$ and $n$, respectively. Let $AB$ denote the sequence $(a_1, \dots, a_m, b_1, \dots, b_n)$ of length $m+n$.

\begin{theorem}\label{thm-join-regular-On2} Let $G$ be a $2r$-regular graph of order $p$ such that $r\ge 1$, $p\ge 3$ and $p-2r\ge 2$ with $\chi(G)=\chi_{la}(G)$. If $G_0=G$ and $G_{k+1} = G_k \vee O_{p-2r}$ for $k\ge 0$, then $G_{k+1}$ is a $(2r+(k+1)(p-2r))$-regular graph of order $p+(k+1)(p-2r)$ with $\chi_{la}(G_k) = \chi_{la}(G) + k$. \end{theorem}

\begin{proof} By Theorem~\ref{thm-join-regular-On}, we only need to consider odd $p$. Let $V(G)=\{u_i\;|\;1\le i\le p\}$ and $V_k(O_{p-2r})=\{v_{k,j}\;|\;1\le j\le p-2r\}$ so that $V(G_k)=V(G_{k-1})\cup V_k(O_{p-2r})$.

\ms\nt Let $\chi(G)=t$ and $f_0$ be a local antimagic $t$-coloring of $G$. Without loss of generality, we may assume $f^+_0(u_1) < f^+_0(u_2) < \cdots < f^+_0(u_t)$ and that for each $a\in[t+1,p]$ there is a $z\in[1,t]$ such that $f^+_0(u_a) = f^+_0(u_z)$. Observe that $f^+_0(u_t)\le q+(q-1)+\cdots+(q-2r+1)=r(2q-2r+1)$, where $q=rp$ is the size of $G_0$.

\ms\nt For convenience, let $p_0=p$ and $q_0=q$, the size of $G$. Clearly, $G_k$ is a $p_{k-1}$-regular graph with order $p_k = p+k(p-2r)$ and size $q_k$, $k\ge 1$. Consider the sequences $U=(u_1,u_2,\ldots, u_p)$ and $V_k = (v_{k,1},v_{k,2},\ldots,v_{k,p-2r})$. In $G_k$, $k\ge 2$, we shall arrange the vertices of $G_{k-1}$ in the sequence $UV_1V_2\cdots V_{k-1}$.

\ms\nt For $k\ge 1$, let $\mathscr P(k)$ be the statement:\\ ``$G_k$ admits a local antimagic $(t+k)$-coloring $f_k$ such that
\begin{enumerate}[(a)]
\item for each $a\in[t+1,p]$, there exists a $z\in[1,t]$ such that $f^+_k(u_a) = f^+_k(u_z)$, and that $f^+_k(v_{i,j})$ dependents on $i$ only;

\item $f^+_{k}(u_1) < f^+_{k}(u_2) <\cdots < f^+_{k}(u_t) < f^+_{k}(v_{1,j}) < \cdots < f^+_{k}(v_{k-1,j}) < f^+_{k}(v_{k,j})$ for $1\le j\le p-2r$."
\end{enumerate}

\nt In $G_1=G_0\vee O_{p-2r}$, we define $f_1: E(G_1)\to [1,q+p(p-2r)]$ such that $f_1(e) = f_0(e)$ if $e\in E(G_0)$ and $f_1(u_iv_{1,j})=m_{i,j}$, where $m_{i,j}$ is the $(i,j)$-entry of a $p\times (p-2r)$ magic rectangle using integers in $[q+1,q+p(p-2r)]$ with constant row sum $R$ and constant column sum $C$. By Lemma~\ref{lem-col-row} we have

\centerline{$C-R= r[2q+1+p(p-2r)]>r(2q-2r+1)\ge f^+_0(u_t).$}

\nt Therefore, $f^+_1(u_i) = f^+_0(u_i)+R$ for $1\le i\le p$, and $f^+_1(v_{1,j})=C$ for $1\le j\le p-2r$. Thus, $f^+_1(u_i)=f^+_1(u_{i'})$ if and only if $f^+_0(u_i) = f^+_0(u_{i'})$ for $1\le i < i'\le p$. Now,  $f^+_1(u_t) =f^+_0(u_t)+R<C= f^+_1(v_{1,j})$. Therefore, $f_1$ is a local antimagic $(t+1)$-coloring of $G_1$. Thus, the statement $\mathscr P(1)$ holds.

\nt We assume that the statement $\mathscr P(s-1)$ holds for $s\ge 2$.

\ms\nt Now, consider $G_s = G_{s-1}\vee O_{p-2r}$.
\begin{enumerate}[(A)]
\item Suppose $s=2n\ge 2$. Note that $G_{2n-1}$ is of order $p_{2n-1}=p+(2n-1)(p-2r)$ which is even. Let $M=(m_{h,j})$ be a $p\times (p-2r)$ magic rectangle using integers in $[q_{2n-1}+1,q_{2n-1}+p(p-2r)]$ and $N=(n_{i,j})$ be a $(p_{2n-1}-p)\times (p-2r)$ magic rectangle using integers in $[q_{2n-1}+p(p-2r)+1,q_{2n-1}+p(p-2r)+(p-2r)(p_{2n-1}-p)]$. Let $c_M$ be the column sum of $M$, and $r_N$ and $c_N$ be then row and column sums of $N$, respectively. So $c_M=\frac{1}{2}p[2q_{2n-1}+p(p-2r)+1]$ and
\begin{align*}
r_N-r_M &=\frac{1}{2}[2q_{2n-1}+1+2p(p-2r)+(p-2r)(p_{2n-1}-p)](p-2r)\\ &\quad -
\frac{1}{2}[2q_{2n-1}+1+p(p-2r)](p-2r).
\end{align*}
Clearly $r_N>r_M$.
By Lemma~\ref{lem-col-row}, we have
\begin{align*}
    c_N-r_N & = \frac{1}{2}[2q_{2n-1}+1+(p-2r)(p+p_{2n-1})](p_{2n-1}-2p+2r)\\ &=\frac{1}{2}[2q_{2n-1}+1+(p-2r)(p+p_{2n-1})](p_{2n-2}-p)
\end{align*}

Let $w_i$ be the $i$-th vertex in $V_1V_2\cdots V_{2n-1}$, $1\le i\le (2n-1)(p-2r)$. We define a bijection $f_{2n}: E(G_{2n-1}\vee O_{p-2r})\to [1,q_{2n-1}+p_{2n-1}(p-2r)]$ by $f_{2n}(e) = f_{2n-1}(e)$ for $e\in E(G_{2n-1})$, $f_{2n}(u_hv_{2n,j})=m_{h,j}$ and $f_{2n}(w_iv_{2n,j}) = n_{i,j}$, $1\le h\le p$, $1\le i\le (2n-1)(p-2r)$, $1\le j\le p-2r$.

Clearly, $f^+_{2n}(u_h) = f^+_{2n-1}(u_h) + r_M$ for $1\le h\le  p$, $f^+_{2n}(v_{i,j}) = f^+_{2n-1}(v_{i,j})+r_N$ for $1\le i\le 2n-1$. Thus the order of $f^+_{s-1}$-values for vertices of $G_{s-1}$ is the same as those of $f^+_{s}$.

Now \begin{align*}& \quad \ 2[f^+_{2n}(v_{2n,j})-f^+_{2n}(v_{2n-1,j})]\\& = 2(c_M + c_N)-2(f^+_{2n-1}(v_{2n-1,j})+r_N)\\ & =2(c_N-r_N)+2c_M-2f^+_{2n-1}(v_{2n-1,j})\\
&\ge [2q_{2n-1}+(p-2r)(p+p_{2n-1})+1](p_{2n-2}-p)\\& \quad +p[2q_{2n-1} +p(p-2r)+1]-(2q_{2n-1}-p_{2n-2}+1)p_{2n-2}\\
& = (p-2r)(p+p_{2n-1})(p_{2n-2}-p)+p^2(p-2r)+p_{2n-2}^2>0.
\end{align*}
So $f^+_{2n}(v_{2n,j}) > f^+_{2n}(v_{2n-1,j})$.

By induction hypothesis and $r_N > r_M$, statement $\mathscr P(2n)$ holds.

\item Suppose $s=2n+1\ge 3$. Note that $G_{2n}$ is of order $p_{2n} = p + 2n(p-2r)$ which is odd. Let $M=(m_{i,j})$ be a $p_{2n}\times (p-2r)$ magic rectangle using integers in $[q_{2n}+1,q_{2n}+p_{2n}(p-2r)]$ with row sum $r_M$ and column sum $c_M$.

Let $z_i$ be the $i$-th vertex in $UV_1V_2\cdots V_{2n}, 1\le i\le p_{2n}$. Define a bijection $f_{2n+1} : E(G_{2n}\vee O_{p-2r})\to [1,q_{2n}+p_{2n}(p-2r)]$ by $f_{2n+1}(e) = f_{2n}(e)$ for $e\in E(G_{2n})$ and $f_{2n+1}(z_iv_{2n+1,j}) = m_{i,j}, 1\le i\le p_{2n}, 1\le j\le p-r$.

Now $f^+_{2n+1}(u_h) = f^+_{2n}(u_h) + r_M$ for $1\le h\le p$, $f^+_{2n+1}(v_{i,j})=f^+_{2n}(v_{i,j})+r_M$ for $1\le i\le 2n$ and $f^+_{2n+1}(v_{2n+1,j}) = c_M$.

Thus the order of $f^+_{s-1}$-values for vertices of $G_{s-1}$ is the same as those of $f^+_{s}$. By a similar calculation as the previous case, we have $f^+_{2n+1}(v_{2n,j}) < f^+_{2n+1}(v_{2n+1,j})$. By induction hypothesis, statement $\mathscr P(2n+1)$ holds.
\end{enumerate}

\nt By mathematical induction, $\mathscr P(k)$ holds for each $k \ge 1$. Thus, $G_k$ admits a local $(t+k)$-coloring $f_k$ with $f^+_{k}(u_1) < f^+_k(u_2) < \cdots < f^+_k(u_t) < f^+_k(v_{1,j}) < \cdots < f^+_k(v_{2k-1,j}) < f^+_k(v_{k,j})$ and that for each $a\in[t+1,p]$, there exists a $z\in[1,t]$ such that $f^+_k(u_a) = f^+_k(u_z)$. Hence, $\chi_{la}(G_k)\ge \chi(G_k) = t+k$, the theorem holds.
\end{proof}

\begin{corollary}\label{cor-GVOn} Let $G$ be a $2r$-regular graph of order $p\ge 5$ in

\centerline{$\{C_n\,|\,n$ odd$\}\cup\{K_5, Q,  K_7, C_7(1,2), \overline{C_3+C_4}, B_2, B_3, B_4, B_5,$ $B_6, C_8(1,2,3)\}$.}

\nt If $G_0=G$ and $G_{k+1} = G_k\vee O_{p-2r}$ for $k\ge 0$, then $G_{k+1}$ is a $(2r+(k+1)(p-2r))$-regular graph of order $p+(k+1)(p-2r)$ with $\chi_{la}(G_k)=\chi_{la}(G)+k$.   \end{corollary}

\begin{theorem}\label{thm-GVCeven} Let $G$ be a $2r$-regular graph of order $p$ such that $r\ge 1, p\ge 3$ and $p-2r+2\ge 2$ with $\chi(G)=\chi_{la}(G)$. If $G_0=G$ and $G_{k+1} = G_k \vee C_{p-2r+2}$ for $k\ge 0$, then $G_{k+1}$ is a $(2r+(k+1)(p-2r+2))$-regular graph of order $p+(k+1)(p-2r+2)$ with $$\chi_{la}(G_k) = \begin{cases}\chi_{la}(G)+2k & \mbox{ if } p \mbox{ is even},\\ \chi_{la}(G) + 3k & \mbox{ otherwise.}\end{cases} $$ \end{theorem}

\begin{proof} Let $V(G)=\{u_i\,|\,1\le i\le p\}$ and $V_k(C_{p-2r+2}) = \{v_{k,j}\,|\,1\le 1\le j\le p-2r+2\}$ so that $E(C_{p-2r+2}) = \{v_{k,j}v_{k,j+1}\,|\,1\le j\le p-2r+2, v_{k,p-2r+3}=v_{k,1}\}$ and $V(G_k) = V(G_{k-1})\cup V_k(C_{p-2r+2})$. Let $\chi(G)=t$ and $f_0$ be a local antimagic $t$-coloring of $G$.  Without loss of generality, we may assume $f^+_0(u_1) > f^+_0(u_2) > \cdots > f^+_0(u_t)$ and that for each $a\in[t+1,p]$ there is a $z\in[1,t]$ such that $f^+_0(u_a) = f^+_0(u_z)$. Note that $q=rp$ is the size of $G_0$. 

\ms\nt For convenience, let $p_0=p$ and $q_0=q$, the size of $G$. Clearly, $G_k$ is a $(p_{k-1}+2)$-regular graph with order $p_k = p+k(p-2r+2)$ and size $q_k$, $k\ge 1$. Consider the sequences $U=(u_1,u_2,\ldots, u_p)$ and $V_k = (v_{k,1},v_{k,2},\ldots,v_{k,p-2r+2})$. In $G_k$, $k\ge 2$, we shall arrange the vertices of $G_{k-1}$ in the sequence $UV_1V_2\cdots V_{k-1}$.

\ms\nt Consider even $p\ge 4$. For $k\ge 1$, let $\mathscr P(k)$ be the statement:\\ ``$G_k$ admits a local antimagic $(t+2k)$-coloring $f_k$ such that
\begin{enumerate}[(a)]
\item for each $a\in[t+1,p]$, there exists a $z\in[1,t]$ such that $f^+_k(u_a) = f^+_k(u_z)$, and that $f^+_k(v_{i,j})$ dependents on $i$ only;
\item $f^+_k(v_{i,2j-1})$ and $f^+_k(v_{i,2j})$ are two distinct constants for $1\le i\le k$ and $1\le j\le (p-2r+2)/2$;
\item $f^+_{k}(u_1) > f^+_{k}(u_2) >\cdots > f^+_{k}(u_t) > f^+_{k}(v_{1,2j}) > f^+_k(v_{1,2j-1}) > \cdots > f^+_{k}(v_{k-1,2j}) > f^+_{k}(v_{k-1,2j-1}) > f^+_{k}(v_{k,2j}) > f^+_{k}(v_{k,2j-1})$ for $1\le j\le (p-2r+2)/2$."
\end{enumerate}

\nt Let $g$ be the local antimagic labeling of $O_{2n}\vee C_{2m}$ as in the proof of Theorem~3.3 in~\cite{LNS-dmgt}. In $G_1=G\vee C_{p-2r+2}$, define a bijection $f_1: E(G_1)\to [1,q+(p+1)(p-2r+2)]$ such that $f_1(e) = g(e)$ if $e\in E(O_p\vee C_{p-2r+2})$, and $f_1(e) = f_0(e) + (p+1)(p-2r+2)$ if $e\in E(G)$. By direct computation, for $1\le j\le (p-2r+2)/2$ and $2m=p-2r+2, 2n=p$, we have $f^+_1(v_{1,2j}) = 4mn^2+12mn-6m+3 > f^+_1(v_{1,2j-1}) = 4mn^2-4mn+2n+10m-1$ and

\begin{align}
f^+_1(u_i) & = f^+_0(u_i) + 2r(p+1)(p-2r+2) + 4m^2n + 4m^2 + m\label{eq-G0VC}\\
 & = f^+_0(u_i) + (2n - 2m + 2)(2n+1)(2m) + 4m^2n + 4m^2 + m\nonumber\\
 & = f^+_0(u_i) + 8mn^2-4m^2n+12mn+5m.\nonumber
\end{align}

\nt Thus, $f^+_1(u_i) = f^+_1(u_{i'})$ if and only if $f^+_0(u_i) = f^+_0(u_{i'})$ for $1\le i < i' \le p$. Since $n>m\ge 2$, we have $f^+_1(u_i) - f^+_1(v_{1,2j}) = f^+_0(u_i) + 4mn(n-m)+ 16mn - 5m - 2n + 1 >0$. Therefore, $f_1$ is a local antimagic $(t+2)$-coloring of $G_1$. Thus, statement $\mathscr P(1)$ holds. We assume that the statment $\mathscr P(s-1)$ holds for $s\ge 2$.  

\ms\nt Now, consider $G_{s} = G_{s-1}\vee C_{p-2r+2}$. Note that $G_{s-1}$ is of order $p_{s-1}=p+(s-1)(p-2r+2)$ which is even. Note that $G_{s-1}$ has regularity $2r+(s-1)(p-2r+2)=2r+p_{s-1}-p$. Define a bijection $f_{s} : E(G_{s-1}\vee C_{p-2r+2}) \to [1, q_{s-1} + (p_{s-1}+1)(p-2r+2)]$ by $f_s(e) = g(e)$ if $e\in E(O_{p_{s-1}}\vee C_{p-2r+2})$, and $f_s(e) = f_{s-1}(e) + (p_{s-1}+1)(p-2r+2)$ if $e\in E(G_{s-1})$. By direct computation, for $1\le j\le (p-2r+2)/2$ and $2m=p-2r+2, 2n=p_{s-1}$, we have $f^+_s(v_{s,2j}) = 4mn^2+12mn-6m+3 > f^+_s(v_{s,2j-1}) = 4mn^2-4mn+2n+10m-1$ and for each $w\in V(G_{s-1})$, by the same computation as \eqref{eq-G0VC} we have

\begin{align*}
f^+_s(w) & = f^+_{s-1}(w) + [2r+p_{s-1}-p](p_{s-1}+1)(p-2r+2) + 4m^2n + 4m^2 + m\\
 & = f^+_{s-1}(w) + 8mn^2-4m^2n+12mn+5m.
\end{align*}

\nt Thus the order of $f^+_{s}$-values for vertices of $G_{s-1}$ is the same as those of $f^+_{s-1}$. Moreover, $f^+_s(v_{s-1,2j-1}) = f^+_{s-1}(v_{s-1,2j-1}) + 8mn^2-4m^2n+12mn+5m > f^+_s(v_{s,2j}) = 4mn^2+12mn-6m+3$.

\nt  Thus, $f^+_s(v_{s-1,2j-1}) >  f^+_s(v_{s,2j})  > f^+_s(v_{s,2j-1})$. Therefore, $f_s$ is a local antimagic $(t+2s)$-coloring of $G_1$. By induction hypothesis, statement $\mathscr P(s)$ holds. 

\nt By mathematical induction, $\mathscr P(k)$ holds for each $k\ge 1$. Thus, $G_k$ admits a local $(t+2k)$-coloring $f_k$. Namely, $f^+_k(v_{i,2j-1})$ and $f^+_k(v_{i,2j'})$ are two distinct constants with $f^+_{k}(u_1) > f^+_{k}(u_2) >\cdots > f^+_{k}(u_t) > f^+_{k}(v_{1,2j}) > f^+_k(v_{1,2j-1}) > \cdots > f^+_{k}(v_{k-1,2j}) > f^+_{k}(v_{k-1,2j-1}) > f^+_{k}(v_{k,2j}) > f^+_{k}(v_{k,2j-1})$ for $1\le j\le (p-2r+2)/2$ and that for each $a\in [t+1,p]$, there exists a $z\in [1,t]$ such that $f^+_k(u_a) = f^+_k(u_z)$. Since $\chi_{la}(G_k)\ge \chi(G_k) = t+2k$, $\chi_{la}(G_k)=t+2k$.

\ms\nt We now consider odd $p\ge 3$. For $k\ge 1$, let $\mathscr P(k)$ be the statement:\\ ``$G_k$ admits a local antimagic $(t+k)$-coloring $f_k$ such that
\begin{enumerate}[(a)]
\item for each $a\in[t+1,p]$, there exists a $z\in[1,t]$ such that $f^+_k(u_a) = f^+_k(u_z)$, and that $f^+_k(v_{i,j})$ dependents on $i$ only;
\item $f^+_k(v_{i,1})$, $f^+_k(v_{i,2j-1})$ and $f^+_k(v_{i,2j'})$ are three distinct constants for $1\le i\le k$, $2\le j\le (p-2r+3)/2$ and $1\le j'\le (p-2r+1)/2$;
\item $f^+_{k}(u_1) < f^+_{k}(u_2) <\cdots < f^+_{k}(u_t) < f^+_{k}(v_{1,1}) < f^+_k(v_{1,2j-1}) < f^+_k(v_{1,2j'}) < \cdots < f^+_{k}(v_{k-1,1}) < f^+_k(v_{k-1,2j-1}) < f^+_k(v_{k-1,2j'}) < f^+_{k}(v_{k,1}) < f^+_k(v_{k,2j-1}) < f^+_k(v_{k,2j'})$ for $1\le i\le k$, $2\le j\le (p-2r+3)/2$ and $1\le j'\le (p-2r+1)/2$."
\end{enumerate}

\ms\nt Let $g$ be the local antimagic $3$-coloring of $C_{p-2r+2}=v_1v_2v_3\cdots v_{p-2r+2}v_1$ provided in \cite{Arumugam}, namely, $g(v_{2j}v_{2j+1})=p-2r+3-j$ and $g(v_{2j-1}v_{2i})=j$ for $1\le j\le (p-2r+3)/2$ where $v_{p-2r+3}=v_1$. Thus, $g^+(v_1)=(p-2r+5)/2 < g^+(v_{2j'})=p-2r+3 < g^+(v_{2j-1})=p-2r+4$ for $1\le j'\le (p-2r+1)/2$ and $2\le j\le (p-2r+3)/2$. In $G_1 = G \vee C_{p-2r+2}$, define $f_1 : E(G_1) \to [1,q+(p+1)(p-2r+2)]$ by
\[f_1(e)=\begin{cases}
f_0(e) &  \mbox{ if }e\in E(G_0),\\
m_{i,j} & \mbox { if } e=u_iv_{1,j},\\
g(e)+q+p(p-2r+2) & \mbox{ if }e\in E(C_{p-2r+2}),
\end{cases}
\]
where $m_{i,j}$ is the $(i,j)$-entry of the $p\times (p-2r+2)$ magic rectangle using integers in $[q+1,q+p(p-2r+2)]$ with constant row sum $R$ and constant column sum $C$.

\ms\nt Now, $f^+_1(u_i) = f^+_0(u_i)+R$ for $1\le i\le p$, and $f^+_1(v_{1,1}) = C + g^+(v_{1,1}) + 2[q+p(p-2r+2)]$ $<$ $f^+_1(v_{2j'}) = C + g^+(v_{1,2j'}) + 2[q+p(p-2r+2)]$ $<$ $f^+_1(v_{1,2j-1}) = C + g^+(v_{1,2j-1}) + 2[q+p(p-2r+2)]$ for $2\le j\le (p-2r+3)/2$ and $1\le j'\le (p-2r+1)/2$. Note that $f^+_0(u_i)\le r(2q-2r+1)$. Moreover, $f^+_1(u_i) = f^+_1(u_{i'})$ if and only if $f^+_0(u_i) = f^+_0(u_{i'})$ for $1\le i < i'\le p$. By Lemma~\ref{lem-col-row}, we have
\begin{align*}
& f^+_1(v_{1,1}) - f^+_1(u_i) \\
 &= C-R - f^+_0(u_i) + 2[q+p(p-2r+2)] + g^+(v_{1,1})\\
 & \ge (r-1)[2q+p(p-2r+2)+1] - r(2q-2r+1) + 2[q+p(p-2r+2)]+ g^+(v_{1,1})\\
  & =rp(p-2r)+2r^2+p^2+2p-1+g^+(v_{1,1})>0.
\end{align*}

\nt Therefore, $f^+_1(v_{1,1}) > f^+_1(u_i)$, and $f_1$ is a local antimagic $(t+3)$-coloring of $G_1$. Thus, statement $\mathscr P(1)$ holds.

\ms\nt We assume the statement $\mathscr P(s-1)$ holds for $s\ge 2$. Now, consider $G_s = G_{s-1}\vee C_{p-2r+2}$.

\ms
\begin{enumerate}[(A)]
  \item Suppose $s=2n\ge 2$. Note that $G_{2n-1}$ is of even order. Let $M=(m_{h,j})$ be a $p\times (p-2r+2)$ magic rectangle using integers in $[q_{2n-1}+1,q_{2n-1}+p(p-2r+2)]$ and $N = (n_{i,j})$ be a $(p_{2n-1}-p)\times (p-2r+2)$ magic rectangle using integers in $[q_{2n-1}+p(p-2r+2)+1,q_{2n-1}+p(p-2r+2)+(p-2r+2)(p_{2n-1}-p)]$. Let $c_M$ be the column sum of $M$, and $r_N$ and $c_N$ be the row and column sums of $N$, respectively. So, $c_M = \frac{1}{2}p[2q_{2n-1}+p(p-2r+2)+1]$ and
\begin{align*}
r_N - r_M &= \frac{1}{2}[2q_{2n-1}+1+2p(p-2r+2)+(p-2r+2)(p_{2n-1}-p)](p-2r+2)\\
 & \quad - \frac{1}{2}[2q_{2n-1}+1+p(p-2r+2)](p-2r+2)
\end{align*}
Clerly $r_N > r_M$. By Lemma~\ref{lem-col-row}, we have
\begin{align*}
c_N - r_N & = \frac{1}{2}[2q_{2n-1}+1+(p-2r+2)(p+p_{2n-1})](p_{2n-1}-2p+2r-2)\\
 & = \frac{1}{2}[2q_{2n-1}+1+(p-2r+2)(p+p_{2n-1})](p_{2n-2}-p)
\end{align*}
Let $w_i$ be the $i$-th vertex in $V_1V_2\cdots V_{2n-1}$, $1\le i\le (2n-1)(p-2r+2)$. We define a bijection $f_{2n} : E(G_{2n-1} \vee C_{p-2r+2})\to [1,q_{2n-1}+(p_{2n-1}+1)(p-2r+2)]$ by $f_{2n}(e) = f_{2n-1}(e)$ for $e\in E(G_{2n-1})$, $f_{2n}(u_hv_{2n,j})=m_{h,j}$ and $f_{2n}(w_iv_{2n,j}) = n_{i,j}$, $1\le h\le p$, $1\le i\le (2n-1)(p-2r+2)$, $1\le j\le p-2r+2$, and $f_{2n}(e) = g(e) + q_{2n-1} + p_{2n-1}(p-2r+2)$.

Clearly, $f^+_{2n}(u_h) = f^+_{2n-1}(u_h) + r_M$ for $1\le h\le p$, $f^+_{2n}(v_{i,j}) = f^+_{2n-1}(v_{i,j})+r_N$ for $1\le i\le 2n-1$. Moreover, $f^+_{2n}(v_{2n,1}) = c_M + c_N + g^+(v_{2n,1}) + 2[q_{2n-1} + p_{2n-1}(p-2r+2)]$ $<$ $f^+_{2n}(v_{2n,2j'}) = c_M + c_N + g^+(v_{2n,2j'}) + 2[q_{2n-1} + p_{2n-1}(p-2r+2)]$ $<$ $f^+_{2n}(v_{2n,2j-1}) = c_M + c_N + g^+(v_{2n,2j-1}) + 2[q_{2n-1} + p_{2n-1}(p-2r+2)]$ for $2\le j\le (p-2r+3)/2$ and $1\le j'\le (p-2r+1)/2$. Now
\begin{align*}
& \quad 2[f^+_{2n}(v_{2n,1}) - f^+_{2n}(v_{2n-1,2j-1})] \\
&= 2\big(c_M + c_N + g^+(v_{2n,1}) + 2[q_{2n-1} + p_{2n-1}(p-2r+2)]\big) - 2(f^+_{2n-1}(v_{2n-1,2j-1})+r_N)\\
 & = 2(c_N - r_N) + 2c_M +  2g^+(v_{2n,1}) + 4[q_{2n-1} + p_{2n-1}(p-2r+2)] - 2f^+_{2n-1}(v_{2n-1,2j-1})\\
 & \ge [2q_{2n-1}+1+(p-2r+2)(p+p_{2n-1})](p_{2n-2}-p) + p[2q_{2n-1}+p(p-2r+2)+1] \\
 & \quad +  2g^+(v_{2n,1})+4[q_{2n-1} + p_{2n-1}(p-2r+2)] - (p_{2n-2}+2)(2q_{2n-1}-p_{2n-2}-1)\\
 & \ge [1+(p-2r+2)(p+p_{2n-1})](p_{2n-2}-p) + p[p(p-2r+2)+1] \\
 & \quad +4p_{2n-1}(p-2r+2) + (p_{2n-2}+2)(p_{2n-2}+1)+ 2g^+(v_{2n,1})>0.
\end{align*}
So, $f^+_{2n}(v_{2n,1}) > f^+_{2n}(v_{2n-1,2j-1})$.
By induction hypothesis and $r_N > r_M$, statement $\mathscr P(2n)$ holds.

\item Suppose $s=2n+1\ge 3$. Note that $G_{2n}$ is of odd order. Let $M=(m_{i,j})$ be a $p_{2n}\times (p-2r+2)$ magic rectangle using integers in $[q_{2n}+1,q_{2n}+p_{2n}(p-2r+2)]$ with row sum $r_M$ and column sum $c_M$.

Let $z_i$ be the $i$-th vertex in $UV_1V_2\cdots V_{2n}$, $1\le i\le p_{2n}$. Define a bijection $f_{2n+1} : E(G_{2n}\vee C_{p-2r+2}) \to [1, q_{2n} + (p_{2n}+1)(p-2r+2)]$ by $f_{2n+1}(e) = f_{2n}(e)$ for $e\in E(G_{2n})$, $f_{2n+1}(z_iv_{2n+1,j}) = m_{i,j}$, $1\le i\le p_{2n}, 1\le j\le p-2r+2$, and $f_{2n}(e) = g(e) + q_{2n}+p_{2n}(p-2r+2)$ for $e\in E(C_{p-2r+2})$.

Now $f^+_{2n+1}(u_h) = f^+_{2n}(u_h)+r_M$ for $1\le h\le p$, $f^+_{2n+1}(v_{i,j}) = f^+_{2n}(v_{i,j})+r_M$ and $f^+_{2n+1}(v_{2n+1,j}) = c_M + g^+(v_{2n+1,j}) + 2(q_{2n}+p_{2n}(p-2r+2))$ for $1\le i\le 2n$ and $1\le j\le p-2r+2$. Hence, $f^+_{2n+1}(v_{2n+1,1}) < f^+_{2n+1}(v_{2n+1,2j'})  < f^+_{2n+1}(v_{2n+1,2j-1})$ for $2\le j\le (p-2r+3)/2$ and $1\le j'\le (p-2r+1)/2$.

By a similar calculation as the previous case, we have $f^+_{2n+1}(v_{2n+1,1}) > f^+_{2n+1}(v_{2n,2j-1})$. By induction hypothesis, statement $\mathscr P(2n+1)$ holds.
\end{enumerate}

\nt By mathematical induction, $\mathscr P(k)$ holds for each $k\ge 1$. Thus, $G_k$ admits a local $(t+3k)$-coloring $f_k$. Namely, $f^+_k(v_{i,1})$, $f^+_k(v_{i,2j-1})$ and $f^+_k(v_{i,2j'})$ are three distinct constants with $f^+_{k}(u_1) < f^+_{k}(u_2) <\cdots < f^+_{k}(u_t) < f^+_{k}(v_{1,1}) < f^+_k(v_{1,2j-1}) < f^+_k(v_{1,2j'}) < \cdots < f^+_{k}(v_{k-1,1}) < f^+_k(v_{k-1,2j-1}) < f^+_k(v_{k-1,2j'}) < f^+_{k}(v_{k,1}) < f^+_k(v_{k,2j-1}) < f^+_k(v_{k,2j'})$ for $1\le i\le k$, $2\le j\le (p-2r+3)/2$ and $1\le j'\le (p-2r+1)/2$ and that for each $a\in [t+1,p]$, there exists a $z\in [1,t]$ such that $f^+_k(u_a) = f^+_k(u_z)$. Since $\chi_{la}(G_k)\ge \chi(G_k) = t+3k$, $\chi_{la}(G_k)=t+3k$.

\ms\nt This completes the proof.\end{proof}

\begin{corollary}\label{cor-GVCn} Let $G$ be a $2r$-regular graph of order $p\ge 3$ in

\centerline{$\{C_n\,|\,n$ odd$\}\cup\{K_5, Q, K_7, C_7(1,2), \overline{C_3+C_4}, B_2, B_3, B_4, B_5,$ $B_6, C_8(1,2,3)\}$.}

\nt If $G_0=G$ and $G_{k+1} = G_k\vee C_{p-2r+2}$ for $k\ge 0$, then $G_{k+1}$ is a $(2r+(k+1)(p-2r+2))$-regular graph of order $p+(k+1)(p-2r+2)$ with $$\chi_{la}(G_k)=\begin{cases}\chi_{la}(G)+2k & \mbox{ if $p$ is even,} \\ \chi_{la}(G)+3k & \mbox{ otherwise}.\end{cases}$$   \end{corollary}

\nt In~\cite{Chai+S+S}, the authors generalized the concept of magic rectangle to nearly magic rectangle as follows. 

\begin{definition} Let $p$ be even and $q$ be odd. A $p\times q$ {\it nearly magic rectangle} is a $p\times q$ matrix that contains each of the integers from the set $[1,pq]$ exactly once, that has constant column sums, and that has row sums that differ by no more than 1. \end{definition}

\begin{theorem}\label{thm-nearmagic} There exists a nearly magic rectangle for all even $p\ge 2$ and odd $q\ge 3$ with constant column sum $p(pq+1)/2$ whereas half of the row sum is $(q(pq+1)-1)/2$ and the other half of the row sum is $(q(pq+1)+1)/2$.   \end{theorem}

\begin{rem}\label{rem-nearmagic} Suppose $N$ is a $p\times q$ nearly magic rectangle and $A$ is a constant. If we add each entry of $N$ by $A$, then the resulting matrix  has  constant column sum $p(pq+1)/2+pA$ whereas half of the row sum is $(q(pq+1)-1)/2+qA$ and the other half of the row sum is $(q(pq+1)+1)/2+qA$.
\end{rem}

\begin{theorem}\label{thm-joinodd} Suppose $m,n\ge 2$ and $t\ge 3$. Let $G$ be a graph of order $2m$ that admits a local antimagic $t$-coloring $f$. If
\begin{enumerate}[(i)]
\item there is a $2$-partition $\{V_1, V_2\}$ of $t$ color classes with $|V_1|=|V_2|=m$,
\item there is no $a,b,a+1,b+1$ such that $a,b+1$ are induced colors of vertices in $V_1$ while $a+1,b$ are induced colors of vertices in $V_2$ and
\item $2n-1> 2m$ or else $m\ge n^2-\frac{3(n-1)}{2}$,
\end{enumerate}
then $\chi_{la}(G\vee O_{2n-1}) \le t+1$. The equality holds if $\chi(G)=t$. \end{theorem}

\begin{proof} Let $V(G)=\{u_i\,|\,1\le i\le 2m\}$ and $V(O_{2n-1})=\{v_j\,|\,1\le j\le 2n-1\}$. By Condition (i), without loss of generality, we can assume $V_1=\{u_i\,|\,1\le i\le m\}$ and $V_2=\{u_i\,|\,m+1\le i\le 2m\}$.
Let $a_{i,j}$ be the $(i,j)$-entry of a $2m\times (2n-1)$ nearly magic rectangle. Without loss of generality, assume row $i$ has row sum $[(2n-1)(2m(2n-1)+1)-1]/2$ when $1\le i\le m$ and  $[(2n-1)(2m(2n-1)+1)+1]/2$ when $m+1\le i\le 2m$. Define a bijection $g:E(G\vee O_{2n-1})\to [1,q+2m(2n-1)]$, where $q=|E(G)|$ such that $g(e)=f(e)$ for $e\in E(G)$ and $g(u_iv_j)=a_{i,j}+q$.

\ms\nt We now have
\begin{enumerate}[(a)]
  \item $g^+(v_j) = 2mq + m(2m(2n-1)+1)$ for $1\le j\le 2n-1$,
  \item $g^+(u_{i})=f^+(u_{i}) + (2n-1)q + [(2n-1)(2m(2n-1)+1)-1]/2$ for $1\le i\le m$, and
  \item $g^+(u_{k})=f^+(u_{k}) + (2n-1)q + [(2n-1)(2m(2n-1)+1)+1]/2$ for $m+1\le k\le 2m$.
\end{enumerate}

\nt Clearly, $f^+(u)=f^+(v)$ implies that $g^+(u)=g^+(v)$, for any $u,v\in V(G)$. Condition (ii) now implies that $g^+(u_{i}) \ne g^+(u_{k})$ for $1\le i\le m, m+1\le k\le 2m$.

\ms\nt Suppose $2n-1> 2m$. Since $f^+(u_i)\ge 1$, we have $g^+(u_i)> g^+(v_j)$ for $1\le i\le 2m$, $1\le j\le 2n-1$.

\nt Suppose $2m\ge 2n$ and $m\ge n^2-\frac{3(n-1)}{2}$. Obviously $q\le m(2m-1)$ and $f^+(u_i)\le \sum\limits_{i=0}^{2m-2}q-i=(2m-1)(q-m+1)$. Therefore, for $1\le i\le m$,
\begin{align*}
2(g^+(v_j)-g^+(u_i))& = [4mq + 2m(2m(2n-1)+1)]\\ &\quad -[2f^+(u_{i}) + 2(2n-1)q + (2n-1)(2m(2n-1)+1)-1]\\
& = (2m-2n+1)(2m(2n-1)+1)+4mq-2f^+(u_{i}) -2(2n-1)q+1\\
& = (2m-2n+1)(2m(2n-1)+1)+(4m-4n+2)q-2f^+(u_{i})+1\\
 &\ge (2m-2n+1)(2m(2n-1)+1)+(4m-4n+2)q-2(2m-1)(q-m+1)+1\\
&=(2m-2n+1)(2m(2n-1)+1)+(-4n+4)q+2(2m-1)(m-1)+1\\
&\ge (2m-2n+1)(2m(2n-1)+1)+(-4n+4)m(2m-1)+2(2m-1)(m-1)+1\\
& = (2m-2n+1)(2m(2n-1)+1)+[(-4n+4)m+2(m-1)](2m-1)+1\\
& = 8m^2 -8mn^2-10m -2n+12mn +4\ge 2m(4m-4n^2+6n-6)+4 \ge 4.
\end{align*}

\nt Since $g^+(u_k) = g^+(u_i)+1$ and $g^+(v_j)-g^+(u_i)\ge 2$ for $1\le j\le 2n-1$, $1\le i\le m$ and $m+1\le k\le 2m$, we also have $g^+(v_j) > g^+(u_k)$. Thus, $g$ is a local antimagic $(t+1)$-coloring of $G\vee O_{2n-1}$ and $\chi_{la}(G\vee O_{2n-1})\le t+1$. If $\chi(G)=t$, we have $\chi_{la}(G\vee O_{2n-1}) \ge \chi(G\vee O_{2n-1}) = t+1$. The equality holds.
\end{proof}

\begin{theorem}\label{thm-QuarticVOodd} Let $A_1=C_7(1,2)$ and $A_2=\overline{C_3+C_4}$ be the quartic graphs in Section~\ref{sec-quartic}. For $n\ge 2$, $\chi_{la}(A_1 \vee O_{2n-1}) = 5$ and $\chi_{la}(A_2\vee O_{2n-1}) = 4$.    \end{theorem}

\begin{proof} Note that $\chi(A_1)=\chi_{la}(A_1)=4$ and $\chi(A_2)=\chi_{la}(A_2)=3$. If $2n-1=3,5$, by Theorem~\ref{thm-regGVOn}, we have $\chi_{la}(G\vee_{2k-1})=5$. If $2n-1\ge 7$, by Theorem~\ref{thm-GVOn}, we have $\chi_{la}(A_1\vee O_{2n-1})=5$. By the same argument, we also have $\chi_{la}(A_2\vee O_{2n-1}) = 4$.
\end{proof}

\begin{theorem}\label{thm-evenVOk} Let $G\in\{C_8(1,2,4), C_8(1,3,4), C_8(1,2,3)\}$. For $k\ge 3$, $\chi_{la}(G\vee O_k) = 5$. \end{theorem}

\begin{proof} In Section~\ref{sec-quintic}, we have $G=C_8(1,2,4)$ admitting a local antimagic 4-coloring $f$ that induced vertex labels $42,49,56,63$ and $\chi(G)=\chi_{la}(G)=4$. Suppose $k\ge 4$ is even. Since $G$ is 5-regular, if $k=2,4,6$, by Theorem~\ref{thm-regGVOn}, we have $\chi_{la}(G\vee O_{k})=5$. Suppose $k\ge 8$, by Theorem~\ref{thm-GVOn}, we also have $\chi_{la}(G\vee O_{k})=5$.

\ms\nt Suppose $k\ge 3$ is odd. It is easy to check that $f$ satisfies Conditions (i) and (ii) of Theorem~\ref{thm-joinodd}. Suppose $k\ge 9$. By Theorem~\ref{thm-joinodd} we have $\chi_{la}(G\vee O_{k}) = 5$.

\ms\nt Suppose $2m=8 > k = 2n-1 = 3,5,7$, then $n=2,3,4$. \begin{enumerate}[(a)]
\item If $n=2$, then $m\ge n^2-\frac{3(n-1)}{2}$. By Theorem~\ref{thm-joinodd} we have $\chi_{la}(G\vee O_{3}) = 5$.
\item If $n=3$, according to the notation defined in the proof of Theorem~\ref{thm-joinodd} we have $g^+(v_j)=324$, $g^+(u_i)=f^+(u_i)+202$ for $1\le i\le 4$ and $g^+(u_i)=f^+(u_i)+203$ for $5\le i\le 8$. Since $f^+(u_i)\in\{42,49,56,63\}$, $\chi_{la}(G\vee O_{5}) = 5$.
\item If $n=4$, according to the notation defined in the proof of Theorem~\ref{thm-joinodd} we have $g^+(v_j)=388$, $g^+(u_i)=f^+(u_i)+340$ for $1\le i\le 4$ and $g^+(u_i)=f^+(u_i)+341$ for $5\le i\le 8$. Similar to part (b), we have $\chi_{la}(G\vee O_{7}) = 5$.
\end{enumerate}

\nt Note that $G=C_8(1,3,4)$ admits a local antimagic 4-coloring $f$ that induced vertex labels $44,48,51,67$ and $G=C_8(1,2,3)$ admits a local antimagic 4-coloring $f$ that induced vertex labels $57,74,77,92$ with $\chi(G)=\chi_{la}(G)=4$, by a similar argument as above, the theorem holds.
\end{proof}

\nt We next consider the join of two non-empty graphs.

\begin{theorem}\label{thm-join-regular} Let $r_1, r_2\ge 2$ and  $t_1, t_2\ge 3$. Suppose $G_i$ is an $r_i$-regular graph of order $p_i$ that admits a local antimagic $t_i$-coloring $f_i$, where $i=1,2$. Let $r=r_1+p_2=r_2+p_1$ so that $G=G_1\vee G_2$ is an $r$-regular graph. For $u\in V(G_1)$ and $v\in V(G_2)$, if $p_1\equiv p_2 \pmod{2}$ and $f_2^+(v) - f_1^+(u)\ne (p_2 - p_1)(p_1p_2+1)/2 + (r_1-r_2)p_1p_2 - p_1r_1r_2/2$, then $\chi_{la}(G) \le t_1+t_2$. The equality holds if $\chi(G_i)=t_i$.
\end{theorem}

\begin{proof} Let $V(G_i) = \{u_{i,j}\,|\,1\le j\le p_i\}$ and $G_i$ has size $q_i$, where $i=1,2$. Without loss of generality, assume $r_1\le r_2$.
 Since $p_i\ge 3$, a $p_1\times p_2$ magic rectangle $(a_{j,k})$ exists using integers in $[1,p_1p_2]$. This $p_1\times p_2$ magic rectangle has row sum constant $s_r=p_2(p_1p_2+1)/2$ and column sum $s_c=p_1(p_1p_2+1)/2$.

\ms\nt Since $G$ has size $q_1 + q_2 + p_1p_2$, define a bijection $g: E(G) \to [1, q_1+q_2+p_1p_2]$ such that
\begin{enumerate}[(i)]
\item  $g(u_{1,j}u_{2,k}) = a_{j,k}$ for $1\le j\le p_1$ and $1\le k\le p_2$,
\item  $g(e) = f_1(e)+p_1p_2$ for $e\in E(G_1)$, and
\item $g(e) = f_2(e) + p_1p_2 + q_1$ for $e\in E(G_2)$.
\end{enumerate}

\ms\nt The induced vertex labels of $G$ under $g$ are given by
\begin{enumerate}[(1)]
  \item $g^+(u_{1,j}) = f_1^+(u_{1,j}) + p_2(p_1p_2+1)/2 + r_1p_1p_2$,
  \item $g^+(u_{2,j}) = f_2^+(u_{2,j}) + p_1(p_1p_2+1)/2 + r_2(p_1p_2+q_1)$.
\end{enumerate}

\nt Clearly, $g^+(u)=g^+(u')$ if and only if $f_i^+(u) = f_i^+(u')$ for $u,u'\in V(G_i)$. Note that $q_i=p_ir_i/2$. The hypothesis then implies that $g^+(u)\ne g^+(v)$ for $u\in V(G_1)$ and $v\in V(G_2)$. Thus, $g$ is a local antimagic $(t_1+t_2)$-coloring. So $\chi_{la}(G)\le t_1+t_2$.

\ms\nt If $\chi(G_i) = t_i$, then $\chi_{la}(G)\ge \chi(G) = t_1+t_2$. The theorem holds.
\end{proof}

\begin{corollary}\label{cor-regVreg} Let $G_1=G$ be an $r$-regular graph of order $p$ such that $\chi_{la}(G) = \chi(G)=\chi$, where $r\ge 2$. For $k\ge 2$, let $G_k = G_{k-1}\vee G_{k-1}$. Then $G_k$ is an $(r+(2^{k-1}-1)p)$-regular graph of order $2^{k-1}p$ with $\chi_{la}(G_k) =\chi(G_k) = 2^{k-1}\chi$.  \end{corollary}

\begin{proof} We first show $\chi_{la}(G_2)=2\chi$. Let $f$ be a local antimagic for $G_1=G$ such that $c(f)=\chi$.
Since $q=pr/2$, $(r+1)r/2\le f^+(u)\le (2q-r+1)r/2$. For $u, v\in V(G)$, $0\le |f^+(u)-f^+(v)|\le (2q-r+1)r/2-(r+1)r/2=(2q-2r)r/2=(p-2)r^2/2<pr^2/2$. Thus the conditions of Theorem~\ref{thm-join-regular} hold. So $\chi_{la}(G_2)=2\chi$.
Note that $G_2$ is an $(r+p)$-regular graph and $\chi(G)=2\chi$. So we may repeat the above argument for $G_k$ when $k\ge 3$.
\end{proof}

\nt In Section~\ref{sec-regular}, we have determined all the regular graphs $G$ of order at most 8 with $\chi_{la}(G)=\chi(G)$. In~\cite[Theorems 3.1, 3.3, and 3.8]{LNS-dmgt}, the authors also proved that $\chi_{la}(C_n\vee O_m)=\chi(C_n\vee O_m) = k$ where $k=3$ for even $m\ge 2$ and even $n\ge 4$, and $k=4$ for odd $m, n\ge 3$. Moreover,  $\chi_{la}(C_n\vee C_n) = \chi(C_n\vee C_n)=6$ for odd $n\ge 5$. In Corollaries~\ref{cor-GVOn} and~\ref{cor-GVCn}, we have found even regular graphs $G$ with $\chi(G) = \chi_{la}(G) \in \{t+s,t+2s,t+3s\}$ for $t\in\{3,4,5,7\}$ and each $s\ge 1$.  By Corollary~\ref{cor-regVreg}, we immediately have the following theorem.

\begin{theorem} For each $k\ge 1$, there exist non-complete $(r+(2^{k-1}-1)p)$-regular graph $G$ of order $2^{k-1}p$ with $\chi_{la}(G) =\chi(G) = 2^{k-1}\chi$, where
\begin{enumerate}[(i)]
  \item $r=2$, $p\ge 3$ is odd, $\chi=3$,
  \item $r=3,4$, $p=6$, $\chi=3$,
  \item $r=4$, $p=7,8$, $\chi=3,4$,
  \item $r=5,6$, $p=8$, $\chi=4$,
  \item $r=n$, $p = 2n-2$, $\chi=3,4$ for $n\ge 4$,
  \item $r=n+2$, $p=2n$, $\chi=6$ for odd $n\ge 5$,
  \item $r=2a+s(b-2a)$, $p=b+(s+1)(b-2a)$, $\chi=t+s$ for $a\in\{1,2,3\}$, $b\in\{5,6,7,8\}$, $t\in\{3,4,5,7\}$, $s\ge 1$,
  \item $r=2a+s(b-2a+2)$, $p=b+s(b-2a+2)$, $\chi=t+2s$ for $a\in\{1,2,3\}$, $b\in\{6,8\}$, $t\in\{3,4,5,7\}$, $s\ge 1$,
  \item $r=2a+s(b-2a+2)$, $p=b+s(b-2a+2)$, $\chi=t+3s$ for $a\in\{1,2,3\}$, $b\in\{3,5,7\}$, $t\in\{3,4,5,7\}$, $s\ge 1$.
\end{enumerate}\end{theorem} 

\begin{corollary} There exist non-complete regular graph $G$ with arbitrarily large order, regularity, and $\chi(G)_{la}=\chi(G)$. \end{corollary}

\section{Conclusion and Open problem}

We have completely determined the local antimagic chromatic number of all connected regular graphs of order at most 8. Consequently, we also show the existence of regular graphs with arbitrarily large order, regularity and local antimagic chromatic numbers. Since every complete graph $K_n, n\ge 3$ is a regular graph with $\chi_{la}(K_n) = |V(K_n)|$, we end this paper with the following question.

\begin{question} Does there exist non-complete regular graph $G$ with $\chi_{la}(G) = |V(G)|$? \end{question}

\end{document}